\documentclass[12pt]{article}
\usepackage{amsmath}
\usepackage{amsbsy}
\usepackage{amssymb}
\usepackage{latexsym}

\def\bbC{{\mathbb C}}

\def\bbZ{{\mathbb Z}}
\def\bbN{{\mathbb N}}

\def\bbT{{\mathbb T}}

\def\bfe{{\boldsymbol E}}

\def\bfp{{\boldsymbol P}}

\def\calc{{\cal C}}

\def\e{{\cal E}}

\def\g{{\cal G}}
\def\h{{\cal H}}
\def\k{{\cal K}}

\begin{document}
\title{\textbf{Subequivalence Relations and Positive-Definite Functions}}
\author{\textbf{A. Ioana, A.S. Kechris, and T. Tsankov}}
\date{}
\maketitle\thispagestyle{empty}

Consider a {\it standard probability space} $(X,\mu )$, i.e., a
space isomorphic to the unit interval with Lebesgue measure.  We
denote by Aut$(X,\mu )$ the automorphism group of $(X,\mu )$, i.e.,
the group of all Borel automorphisms of $X$ which preserve $\mu$
(where two such automorphisms are identified if they are equal
$\mu$-a.e.).  A Borel equivalence relation $E\subseteq X^2$ is
called {\it countable} if every $E$-class $[x]_E$ is countable and
{\it measure preserving} if every Borel automorphism $T$ of $X$ for
which $T(x)Ex$ is measure preserving.  Equivalently, $E$ is
countable, measure preserving iff it is induced by a measure
preserving action of a countable (discrete) group on $(X,\mu )$ (see
Feldman-Moore \cite{FM}).

To each countable, measure preserving equivalence relation $E$ one
can assign the positive-definite function $\varphi_E(S)$ on
Aut$(X,\mu )$ given by $\varphi_E (S)=\mu (\{x:S(x)Ex\})$; see
Section 1. Intuitively, $\varphi_E(S)$ measures the amount by which
$S$ is ``captured'' by $E$.  This positive-definite function
completely determines $E$.

We use this function to measure the proximity of a pair $E\subseteq
F$ of countable, measure preserving equivalence relations.  In
Section 2, we show, among other things, the next result, where we
use the following notation:  If a countable group $\Gamma$ acts on
$X$, we also write $\gamma$ for the automorphism
$x\mapsto\gamma\cdot x$; if $A\subseteq X$ and $E$ is an equivalence
relation on $X$, then $E|A=E\cap A^2$ is the restriction of $E$ to
$A$; if $E\subseteq F$ are equivalence relations, then $[F:E]=m$
means that every $F$-class contains exactly $m$ classes; if $F$ is a
countable, measure preserving equivalence relation on $(X,\mu )$,
then $[F]$ is the {\it full group} of $F$, i.e., $[F]=\{T\in{\rm
Aut}(X,\mu ):T(x)Fx, \mu {\rm -a.e.} (x)\}$.
\medskip

\noindent{\bf Theorem 1}. {\it Let $\Gamma$ be a countable group and consider a
measure preserving action on $\Gamma$ on $(X,\mu )$ with induced equivalence
relation $F=E^X_\Gamma$.

i) If $E\subseteq F$ is a subequivalence relation and
$\inf_{\gamma\in\Gamma} \varphi_E(\gamma )=\varphi^0_E>0$, then
there is an $E$-invariant Borel set $A\subseteq X$ of positive
measure such that $[F|A:E|A]=m\leq\tfrac{1}{\varphi^0_E}$, so that
if $\varphi^0_E>\tfrac{1}{2}$, $F|A=E|A$.

ii) If $E$ is any countable, measure preserving equivalence relation
and $\epsilon >0$, then $\forall\gamma\in\Gamma (\varphi_E(\gamma
)\geq 1-\epsilon )$ implies $\forall S\in [F] (\varphi_E(S)\geq
1-4\epsilon )$.}

\medskip
{\bf Remark.} Popa pointed out that some version of part (i) of the
preceding theorem was known in the theory of operator algebras, see,
for example, the appendix to Popa \cite{PO1}. Actually our initial
proof of that theorem was inspired by Popa's technique of
conjugating subalgebras in a finite von Neumann algebra (see Section
2 in \cite{PO2}) but, for consistency with the rest of the paper, we
give another self-contained, ergodic-theoretic proof.

\medskip
With some additional work, Theorem 1 has the following consequences:

a) In the context of i), if $\varphi^0_E > 0$, the action of
$\Gamma$ is {\it free} (i.e., $\gamma\cdot x\neq x,\forall\gamma\neq
1$) and $E$ is induced by a free action of a countable group
$\Delta$, then $\Gamma ,\Delta$ are measure equivalent (ME).

b) Again in the context of i), if $\varphi^0_E>\tfrac{1}{2}$ and $E$
is aperiodic (i.e., has no finite classes), then $C_\mu (F)\leq
C_\mu (E)$, where $C_\mu (R)$ is the cost of an equivalence relation
$R$ (see [G1] or [KM] for the theory of costs).

c) In i) if $\varphi^0_E>\tfrac{3}{4}$, then we can find $A$ so that $\mu (A)\geq
4\varphi^0_E-3$.

In Section 3, we consider a recent co-inducing construction of
Epstein \cite{E}. Given a measure preserving, ergodic action $b_0$
of a countable group $\Gamma$ on $(X,\mu )$ with associated
equivalence relation $F=E^X_\Gamma$ and a free, measure preserving
action $a_0$ of a countable group $\Delta$ on $(X,\mu )$ with
associated equivalence relation $E=E^X_\Delta\subseteq
F=E^X_\Gamma$, Epstein's construction gives for any measure
preserving action $a$ of $\Delta$ on a space $(Y,\nu )$, a measure
preserving action $b$ of $\Gamma$ on a space $(Z,\rho )$, called the
{\it co-induced action of a modulo} $(a_0,b_0 )$, in symbols $b={\rm
CInd}(a_0, b_0)^\Gamma_\Delta (a)$. This construction has important
applications in the study of orbit equivalence of actions -- see
Epstein \cite{E}.

For further potential applications of this method, it seems that one
should have a better understanding of the connection of ergodic
properties between $a,b$ as above.  We show, for example, that if
$b_0$ is free, mixing and $a_0$ is ergodic, then: $a$ is mixing
$\Rightarrow b$ is mixing.  There are however interesting situations
under which $b$ is always mixing for arbitrary $a$.  It turns out
that this phenomenon, for given $(a_0,b_0)$, is connected to the
positive-definite function discussed earlier. We show the following:

\medskip
\noindent{\bf Theorem 2}. {\it If $b_0$ is mixing, the following are
equivalent:

(i) For all actions $a$ of $\Delta$, $b={\rm
CInd}(a_0,b_0)^\Gamma_\Delta (a)$ is mixing,

(ii) $\varphi_E(\gamma )\rightarrow 0$ as
$\gamma\rightarrow\infty$.}

\medskip
The condition (ii) in Theorem 2 somehow asserts that $E$ is
``small'' relative to $F$.  In the opposite case we have the
following fact.  If $\varphi^0_E = \inf_{\gamma\in\Gamma}\varphi_E
(\gamma )>0$, then $b$ is ergodic $\Rightarrow a$ is ergodic.

It is well-known that for any ergodic $b_0$ as above one can find a
free, mixing action $a_0$ of $\Delta =\bbZ$ with $E\subseteq F$
(see, e.g., Zimmer \cite{Z}, 9.3.2). We show that when $b_0$ is
mixing, one can find such an $a_0$ so that (ii) of Theorem 2 holds.
This gives a method of producing, starting with arbitrary measure
preserving $\bbZ$ actions, apparently new types of measure
preserving, mixing actions of any infinite group $\Gamma$.

\medskip
\noindent{\bf Theorem 3}. {\it Let $\Gamma$ be an infinite countable
group, and let $b_0$ be a free, measure preserving, mixing action of
$\Gamma$ on $(X,\mu )$.  Then there is a free, measure preserving,
mixing action $a_0$ of $\bbZ$ on $(X,\mu )$ such that
$E=E^X_\bbZ\subseteq F=E^X_\Gamma$ and $\varphi_E(\gamma
)\rightarrow 0$ as $\gamma\rightarrow\infty$.}

\medskip
When the group $\Gamma$ is non-amenable, then by work of
Gaboriau-Lyons \cite{GL} one can find a free, mixing action $b_0$ of
$\Gamma$ on $(X,\mu )$ and a free ergodic action $a_0$ of $F_2$ on
$(X,\mu )$ with $E=E^X_{F_2}\subseteq F=E^X_\Gamma$.  We show again
that such $a_0$ can be found so that (ii) of Theorem 2 holds.  This
is joint work with I. Epstein.

\medskip
\noindent{\bf Theorem 4 (with I. Epstein)}.  {\it Let $\Gamma$ be a
non-amenable countable group.  Then there is a free, measure
preserving, mixing action $b_0$ of $\Gamma$ on $(X,\mu )$ and a
free, measure preserving, ergodic action $a_0$ of $F_2$ on $(X,\mu
)$ such that $E=E^X_{F_2}\subseteq F=E^X_\Gamma$ and
$\varphi_E(\gamma ) \rightarrow 0$ as $\gamma\rightarrow\infty$.}

\medskip
{\bf Remark.} Our proof of Theorem 4 uses the Gaboriau-Lyons
\cite{GL} result and additionally the co-inducing construction to
produce a pair of actions satisfying Theorem 4. In \cite{GL} the
authors actually produce two different pairs of actions as above.
After seeing a preliminary version of our paper, Lyons pointed out
that their first construction can be shown to satisfy Theorem 4,
using results of Benjamini-Lyons-Peres-Schramm \cite{BLPS}, in
particular formula (13.8) in that paper. Subsequently, we realized
that the second construction of \cite{GL} also may give a pair of
equivalence relations satisfying Theorem 4. More precisely, if one
chooses a Cayley graph of $\Gamma$ with sufficiently many generators
and $p \in (p_c, p_u)$ close enough to $p_c$, then the
subequivalence relation one obtains using the method of \cite{GL}
and our Lemma 4.2 will, in fact, satisfy Theorem 4. See \cite{GL},
Pak--Smirnova-Nagnibeda \cite{PS} and the proof of Benjamini-Schramm
\cite{BS}, Theorem 4, for more details. We also note that sometimes
the cluster subequivalence relation $E$ for Bernoulli percolation in
the non-uniqueness phase does not satisfy $\varphi_E(\gamma) \to 0$
as $\gamma \to \infty$ (see Lyons-Schramm \cite{LS}, Remark 1.3).

\medskip
We use Theorem 4 to study the complexity of the classification
problem of free, measure preserving, ergodic actions of a countable
group $\Gamma$ under {\it orbit equivalence} (OE). After a series of
earlier results that dealt with various important classes of
non-amenable groups (see Gaboriau-Popa \cite{GP}, Hjorth \cite{H2},
Ioana \cite{I}, Kida \cite{KI}, Monod-Shalom \cite{MS}, Popa
\cite{PO2}), Epstein \cite{E} finally showed that in general any
non-amenable group admits uncountably many non-orbit equivalent
free, measure preserving, ergodic actions. This was proved earlier
by Ioana \cite{I} in the case where $F_2\leq\Gamma$, and his main
lemma in that proof could be also used to derive, in this case, the
stronger fact that the equivalence relation $E_0$ (or $2^\bbN$,
where $xE_0y\Leftrightarrow\exists n\forall m\geq n (x(m)=y(m)$) can
be Borel reduced to OE on the space of free, measure preserving,
ergodic actions of $\Gamma$. Moreover OE on that space cannot be
classified by countable structures (see \cite{K}, Section 17, {\bf
(B)}). However, it was not known whether this non-classification
result extends to {\it all} non-amenable groups and whether every
non-amenable group admits uncountably many non-orbit equivalent
free, measure preserving, {\it mixing} actions. Putting together
Theorems 2,4 and the work of Epstein \cite{E} leads now to the
following positive answer.  This is again a joint result with I.
Epstein.

\medskip
\noindent{\bf Theorem 5 (with I. Epstein)}.  {\it Let $\Gamma$ be a
non-amenable countable group.  Then $E_0$ can be Borel reduced to
{\rm OE} on the space of free, measure preserving, mixing actions of
$\Gamma$ and {\rm OE} in this space cannot be classified by
countable structures.}

\medskip
Thus we have the following strong dichotomy concerning orbit
equivalence: If $\Gamma$ is (infinite) amenable, there is exactly
one free, measure preserving, ergodic action of $\Gamma$ up to OE,
while if $\Gamma$ is non-amenable, OE of free, measure preserving,
mixing actions of $\Gamma$ is unclassifiable in a very strong sense.

The proof of Theorem 5 shows that the conclusion in that theorem
also holds if OE is replaced by conjugacy (isomorphism) of actions.
This fact is also known to be true for abelian $\Gamma$ (see
\cite{K}, 5.7, where the proof is presented for $\bbZ$ but easily
generalizes to any abelian $\Gamma$).

In Section 4, we review some basic facts concerning invariant bond
percolation on Cayley graphs of finitely generated groups (see
Lyons-Schramm \cite{LS} or Lyons-Peres \cite{LP}). We also give in
Section 5, {\bf (C)} an alternative proof of Theorem 4.1 (for Cayley
graphs) in Lyons-Schramm \cite{LS}, using our Theorem 1.

In Section 5, we apply the preceding results to property (T) groups.
Recall that a {\it Kazhdan pair} $(Q,\epsilon )$ for such a group
consists of a finite generating set $Q\subseteq\Gamma$ and a
positive $\epsilon$ such that for any unitary representation $\pi$
of $\Gamma$ on a Hilbert space $\h$, if there is a vector $\xi\in\h$
with $\Vert\pi (\gamma )(\xi )-\xi\Vert <\epsilon\Vert\xi\Vert ,
\forall\gamma\in Q$, then there is a non-0 invariant vector.  We
state below some sample results.
\medskip

\noindent{\bf Theorem 6}.  {\it Let $\Gamma$ be an infinite group
with property {\rm (T)}, $(Q,\epsilon )$ a Kazhdan pair and $\bfp$
an invariant, ergodic, insertion-tolerant bond percolation on the
Cayley graph $\g_Q$ of $\Gamma$ {\rm (}with respect to $Q${\rm )}.
If the survival probability $\bfp (\{\omega :\omega (e)=1\})$ of
each edge $e$ is $>1-\tfrac{\epsilon^2}{2}$, then $\omega$ has a
unique infinite cluster, $\bfp$-a.s.  In particular, if $p_u(Q)$ is
the critical probability for existence of unique infinite clusters
in Bernoulli percolation on this Cayley graph, then
$p_u(Q)\leq1-\tfrac{\epsilon^2}{2}.$}

\medskip
\noindent{\bf Theorem 7}. {\it For each $\rho >0$ and every infinite
group $\Gamma$ with property {\rm (T)}, there is a finite set of
generators $Q$ for $\Gamma$ such that for any invariant, ergodic,
insertion-tolerant bond percolation $\bfp$ on $\g_Q$, if the
survival probability of each edge is $\geq\rho$, then $\omega$ has a
unique infinite cluster, $\bfp$-a.s.}

\medskip
{\bf Remark.} Lyons-Schramm \cite{LS} had earlier shown that, in the
notation of Theorem 6, $p_u(Q) < 1$. Lyons pointed out that one
could also easily deduce a version of Theorem 6 from the results of
their paper (with perhaps a different constant instead of
$1-\frac{\epsilon^2}{2})$. Similarly for Theorem 5.9 below. Finally,
Lyons mentions that for the special case $\bfp = \bfp_p$, Bernoulli
percolation, Theorem 7 was known even for groups $\Gamma$ for which
there exists $Q$ such that $p_u(Q) < 1$ but which do not necessarily
satisfy property (T).

\medskip
Denote below by $C(\Gamma )$ the cost of a countable group $\Gamma$.
If $\Gamma$ is an infinite countable group with property (T) and $n$
is the smallest cardinality of a set of generators for $\Gamma$,
then we have $1\leq C(\Gamma )<n$. At this time no example of a
property (T) group with $C(\Gamma )>1$ is known. We obtain here some
upper bounds for $C(\Gamma )$ in terms of $n,\epsilon$, where
$n={\rm card} (Q)$ and $(Q,\epsilon )$ is a Kazhdan pair.  One such
result is the following:
\medskip

\noindent{\bf Theorem 8}. {\it Let $\Gamma$ be an infinite group
with property {\rm (T)} and $(Q,\epsilon )$ a Kazhdan pair for
$\Gamma$. If ${\rm card} (Q)=n$, then
\[
C(\Gamma )\leq n\left (1-\frac{\epsilon^2}{2}\right )+
\frac{n-1}{2n-1}.
\]}

Another example is the following.

\medskip
\noindent{\bf Theorem 9}. {\it Let $\Gamma$ be an infinite group
with property {\rm (T)} and let $(Q,\epsilon )$ be a Kazhdan pair,
where $Q$ contains an element of infinite order.  Then if {\rm
card}$(Q)=n$,
\[
C(\Gamma )\leq n-\frac{\epsilon^2}{2}.
\]
In particular, if $\Gamma$ is torsion-free and 2-generated, then $C(\Gamma )\leq 2-
\tfrac{(\epsilon_2)^2}{2}$, where $\epsilon_2$ is the sup of the $\epsilon$ such that
$(Q,\epsilon )$ is a Kazhdan pair with {\rm card}$(Q)=2$.}

\medskip
{\bf Remark.}  Since in this paper we work completely in a measure
theoretic context, we neglect null sets if there is no danger of
confusion.  So given a measure space $(X,\mu )$, we do not often
distinguish between a statement being true for all $x\in X$ or for
all $x\in X,\mu$-a.e.
\medskip

{\bf Acknowledgements}.  The research of A.S.K. and T.T. was
partially supported by NSF Grant DMS-0455285. We would like to thank
I. Epstein for allowing us to include here our joint results. We
would also like to thank R. Lyons, S. Popa and Y. Shalom for many
valuable comments.

\section{Equivalence relations and positive-definite functions on Aut$(X,\mu )$}

Let $(X,\mu)$ be a standard measure space and Aut$(X,\mu )$ the group of measure
preserving automorphisms of $(X,\mu )$.  Denote by $u$ the uniform topology on
Aut$(X,\mu )$, induced by the metric
\[
\delta_u (S,T)=\mu (\{x:S(x)\neq T(x)\}).
\]
Let $E\subseteq X^2$ be a countable, measure preserving equivalence relation on $X$.
Define on Aut$(X,\mu )^2$:
\[
\psi_E(S,T)=\mu (\{x:S^{-1}(x)ET^{-1}(x)\}).
\]
So if $E=\Delta$, the equality relation, then
$1-\psi_E(X,T)=\delta_u(S,T)$. We claim that $\psi_E$ is a
continuous, positive-definite function on $({\rm Aut}(X, \mu ),u)$.
Recall that a function $\psi : Y\times Y \rightarrow \bbC$ is {\it
positive-definite} if for every finite subset $\{ y_1, \dots , y_n\}
\subseteq Y$ and every $\alpha_i \in \bbC, 1\leq i \leq n,$ we have
$\sum_{1\leq i,j\leq n} \bar\alpha_i \alpha_j \psi (y_i, y_j) \geq
0$. The proof is similar to that in Aizenman-Newman \cite{AN}. Fix a
finite set $\{S_1,\dots ,S_n\}\subseteq{\rm Aut}(X,\mu )$ and
$\alpha_i\in\bbC$, $1\leq i\leq n$, in order to show that
\[
\sum_{1\leq i,j\leq n}\bar\alpha_i\alpha_j\psi_E(S_i,S_j)\geq 0.
\]
For each $x\in X$, define the equivalence relation $\sim_x$ on
$\{1,\dots ,n\}$ by
\[
i\sim_xj\Leftrightarrow S^{-1}_i(x)ES^{-1}_j(x).
\]
Let $C^x_1,\dots ,C^x_m$ be the $\sim_x$-classes.  Then if (1) denotes the above
sum, we have
\begin{align*}
(1) &=\int\sum_{1\leq i,j\leq n}\bar\alpha_i\alpha_j\chi_{\{x:S^{-1}_i(x)E
S^{-1}_j(x)\}}d\mu\\
&=\int\sum^m_{k=1}(\sum_{i,j\in C^x_k}\bar\alpha_i\alpha_j)d\mu (x)\\
&=\int\sum^m_{k=1}|\sum_{i\in C^x_k}\alpha_i|^2d\mu (x)\\
&\geq 0.
\end{align*}

Thus $1-\psi_E(S,T)$ is negative-definite. Recall that a function
$\rho : Y\times Y \rightarrow \bbC$ is {\it (conditionally)
negative-definite} if $\rho (y,z) = \overline{\rho (z,y)}$ and for
every finite subset $\{ y_1, \dots , y_n\} \subseteq Y$ and every
$\alpha_i \in \bbC, 1\leq i \leq n,$ with $\sum^n_{i=1} \alpha_i =
0$, we have $\sum_{1\leq i,j\leq n} \bar\alpha_i \alpha_j \rho (y_i,
y_j) \leq 0$. In particular, if $E=\Delta$, then
\[
\delta_u(S,T)=1-\psi_\Delta (S,T),
\]
so the metric $\delta_u$ is negative-definite.

Note also that $\psi_E$ is left-invariant, so
\[
\varphi_E(S)=\psi_E(1,S)
\]
is a continuous, positive-definite function on $({\rm Aut}(X,\mu
),u)$. Recall again that a function $\varphi : G \rightarrow \bbC$
on a group $G$ is {\it positive-definite} if the function $\psi :
G\times  G \rightarrow \bbC$ defined by $\psi (g,h) = \varphi
(g^{-1} h)$ is positive-definite. If $A_E(S)=\{ x:S(x)Ex\}$, then
$\varphi_E(S)=\mu (A_E(S))$, and we view the quantity $\varphi_E(S)$
as measuring the amount by which $S$ is ``captured" by $E$.  By the
GNS construction there is a (unique) triple $(\pi_E,\h_E,\xi_E)$,
consisting of a cyclic continuous representation of $({\rm
Aut}(X,\mu ),u)$ on a Hilbert space $\h_E$ with cyclic unit vector
$\xi_E\in\h_E$ such that
\[
\varphi_E(S)=\langle\pi_E(S)(\xi_E),\xi_E\rangle .
\]
Now if $[E]$ is the full group of $E$, then
\begin{align*}
S\in [E]&\Leftrightarrow\varphi_E(S)=1\\
&\Leftrightarrow\pi_E(S)(\xi_E)=\xi_E,
\end{align*}
i.e., $[E]$ is the stabilizer of $\xi_E$ in $\pi_E$.  In particular,
$\varphi_E$ completely determines $[E]$ and thus $E$, i.e., $E$ is
encoded in $\varphi_E$.

It is not clear how to characterize the continuous, positive-definite functions
$\varphi$ on $({\rm Aut}(X,\mu ),u)$, which are of the form $\varphi_E$ for
some $E$.  Clearly any such $\varphi$ satisfies $0\leq\varphi\leq 1$ and $\varphi
(1)=1$.  Another necessary condition is that ker$(\varphi )=\{S\in{\rm Aut}(X,
\mu ):\varphi (S)=1\}$ (which is a closed subgroup of $({\rm Aut}(X,\mu ),u)$)
is separable in the uniform topology.  The following observation
may also be relevant here.  Let $\Gamma\leq {\rm ker}(\varphi )$ be a
countable dense subgroup
of ker$(\varphi )$.  If $\varphi$ is of the form $\varphi_E$ for some $E$,
then $\Gamma$ is uniformly dense in ${\rm ker}(\varphi )={\rm ker}(\varphi_E)
=[E]$, so $E=E^X_\Gamma =$ the equivalence relation induced by $\Gamma$.

Next consider the negative-definite function
\[
\theta_E(S)=1-\varphi_E(S)
\]
on Aut$(X,\mu )$. Recall that if $\theta : G\rightarrow\bbC$ is a
function on a group $G$, then $\theta$ is negative-definite if the
function $\rho (g,h) = \theta (g^{-1} h)$ is negative-definite.

Put also
\[
\delta_u(S,[E])= \inf \{\delta_u(S,T):T\in [E]\},
\]
for the distance (in $\delta_u$) of $S$ to $[E]$.  Then we have
\medskip

\noindent{\bf Proposition 1.1}. {\it $\theta_E(S)=\delta_u(S,[E])=\inf\{\delta_u
(S,T):T\in [E]\}$ and moreover this $\inf$ is attained.}
\medskip

{\bf Proof}.  We will use the following well-known fact:
\medskip

\noindent{\bf Lemma 1.2}. {\it Let $S\in{\rm Aut}(X,\mu )$ and let $E$ be a
countable, measure preserving equivalence relation on $X$.  Then there is
$T\in [E]$ such that $S(x)=T(x)$, whenever $S(x)Ex$.}
\medskip

Granting this, given any $S$, find $T$ as in 1.2 and note that
\begin{align*}
\delta_u(S,T)&=\mu (\{x:\neg \; S(x) E x\})\\
&=\theta_E(S).
\end{align*}
On the other hand, for any $R\in [E]$, $\{x:\neg\; S(x) E
x\}\subseteq \{x:S(x)\neq R(x)\}$, so
$\theta_E(S)\leq\delta_u(S,R)$, thus
$\theta_E(S)=\delta_u(S,T)=\delta_u (S,[E])$.

\medskip
{\bf Proof of 1.2}.  Let $A=\{x:S(x)Ex\}$ and $B=S(A)$.  It is enough to find a
Borel bijection (modulo null sets) $S':A\cup B\rightarrow A\cup B$ with $S'(x)=
S(x)$ for $x\in A$ and $S'(x)Ex$, for $x\in A\cup B$.  Then we can take $T=S'
\cup {\rm id}|(X\setminus (A\cup B))$.

Put $Y=A\cup B$ and consider the equivalence relation $F$ on $Y$ induced by $S|A$.
Some $F$-classes $C$ will consist of a cycle $\{x,S(x),\dots ,S^n(x)\}$, where
$S^{n+1}(x)=x$. For such $C$, we have $C\subseteq A$, so we can let $S'(x)
=S(x),\forall x\in C$.  In every other $F$-class $C$, we can define the ordering
\[
x<_Cy\Leftrightarrow\exists n>0(S^n(x)=y).
\]
The union of the infinite $C$ in which there is a largest or
smallest element in $<_C$ has clearly measure 0.  So we can assume
that $<_C$ is either a finite ordering, with largest and smallest
elements $b_C, a_C$, resp., in which case $A\cap C=C\setminus
\{b_C\}$ or else $<_C$ looks like a copy of the order on $\bbZ$, in
which case $C\subseteq A$.  In the first case, we define $S'$ on $C$
by $S'(x)=S(x)$, if $x\neq b_C$ and $S'(b_C)=a_C$.  In the second
case, we put $S'(x)=S(x),\forall x\in C$.  This clearly
works.\hfill$\dashv$
\medskip

Note that if $\delta_u$ also denotes the metric induced by $\delta_u$ on the
homogeneous space Aut$(X,\mu )/[E]$, i.e.,
\begin{align*}
\delta_u(S[E],T[E])&=\inf \{\delta_u(S',T'):S'\in S[E], T'\in T[E]\}\\
&=\delta_u(S,T[E])=\delta_u(T,S[E]),
\end{align*}
then
\begin{align*}
\delta_u(S[E],T[E])&=\delta_u(S^{-1}T,[E])\\
&=\theta_E(S^{-1}T)=1-\psi_E(S,T),
\end{align*}
so $\delta_u$ on Aut$(X,\mu )/[E]$, with the quotient topology of
$u$, is a continuous, negative-definite function.

And we conclude with some further observations on metrics on Aut$(X,\mu )$ and
certain subgroups of it.

The weak topology $w$ on Aut$(X,\mu )$ is induced by the metric
\[
\delta_w(S,T)=\sum^\infty_{n=1}2^{-n}\mu (S(A_n)\Delta T(A_n)),
\]
where $\{A_n\}$ is dense in the measure algebra MALG$_\mu$ of $(X,\mu )$.  Now for
each fixed Borel set $A\subseteq X$,
\[
\rho_A(S,T)=\mu (S(A)\Delta T(A))
\]
is negative-definite, since
\begin{align*}
\rho_A(S,T)&=\int |\chi_{S(A)}-\chi_{T(A)}|^2d\mu\\
&=\Vert\chi_{S(A)}-\chi_{T(A)}\Vert^2_2,
\end{align*}
and the function $(\xi ,\eta )\mapsto\Vert\xi -\eta\Vert^2_2$ is
negative-definite on $L^2(X,\mu )$.  It follows that the left-invariant metric
$\delta_w$ is negative-definite.  In particular, the complete metric $\bar\delta_w
(S,T)=\delta_w(S,T)+\delta_w(S^{-1},T^{-1})$ on Aut$(X,\mu )$ is also
negative-definite.

Now consider an aperiodic (i.e., having infinite classes) $E$ and the normalizer
$N[E]$ of its full group.  Then $N[E]$ has a canonical topology induced by the
complete metric
\[
\bar\delta_{N[E]}(S,T)=\bar\delta_w(S,T)+\sum^\infty_{n=1}2^{-n}\delta_u
(S\gamma_nS^{-1},T\gamma_nT^{-1}),
\]
where $\{\gamma_n\}$ is a countable subgroup of Aut$(X,\mu )$ inducing $E$
(see, e.g., Kechris \cite{K}).  Since for each $n$, the function $(S,T)\mapsto\delta_u
(S\gamma_nS^{-1},T\gamma_nT^{-1})$ is negative-definite, so is $\bar\delta_{N[E]}
(S,T)$ on $N[E]$.

\section{Proximity of subequivalence relations}

{\bf (A)} We view the quantity $\varphi_E(S)=\mu (\{x:S(x)Ex\})$ as measuring the
amount by which $S$ is captured by $E$.  We will next see that if a countable group
$\Gamma$ acts in a measure preserving way on $(X,\mu )$ inducing an equivalence
relation $F=E^X_\Gamma$, and every element of $\Gamma$ (viewed as an element of
Aut$(X,\mu )$ via $x\mapsto\gamma\cdot x,\gamma\in\Gamma$) is ``substantially
captured'' by $E$, then $E,F$ are somehow ``close'' to each other.

Towards this goal we will study a canonical representation
associated to a pair $E\subseteq F$ of countable measure preserving
equivalence relations on $(X,\mu )$.

Let such $E,F$ be given and decompose $X=\bigsqcup_{N\in\{1,2,\dots
,\aleph_0\}} X_N$, where $X_N=\{x:$ There are exactly $N\ E$-classes
in $[x]_F\}$, so that $X_N$ is $F$-invariant.  Thus
$[F|X_N:E|X_N]=N$.  If $F$ is ergodic, clearly $X=X_N$ for some $N$.
Fix now for each $N$ a sequence of Borel functions $\{C^{(N)}_n
\}_{n\in N}, C^{(N)}_n:X_N\rightarrow X_N$, where we identify $N$
here with $\{0,\dots ,N-1\}$, if $N$ is finite, and with $\bbN$, if
$N=\aleph_0$, such that $C^{(N)}_0= {\rm id}|X_N$, for each $x\in
X_N, C^{(N)}_n(x)\neq C^{(N)}_m(x)$, if $m\neq n$, and
$\{C^{(N)}_n(x)\}$ is a transversal for the $E$-classes contained in
$[x]_F$.  These are called {\it choice functions}.
\medskip

{\bf Remark 2.1}. For further reference, notice that if $E$ is
ergodic, so that $X=X_N$ for some $N$, then we can take the choice
functions $C^{(N)}_n=C_n$ to be 1-1, i.e., to be in Aut$(X,\mu )$.
To see this, start with arbitrary $\{C^{(N)}_n\}= \{C_n\}$.  Fix
$n\in N$ and consider $C_n$.  As it is countable-to-1, let $X=
\bigsqcup^\infty_{k=1}Y_k$ be a Borel partition such that $C_n|Y_k$
is 1-1.  Let then $Z_k=C_n(Y_k)$, so that $\mu (Z_k)=\mu (Y_k)$.
Since $E$ is ergodic, there is $T_k\in [E]$ with $T_k(Z_k)=Y_k$. Let
then $D_n(x)=T_k (C_n(x))$, if $x\in Y_k$. We have
$D_n(x)EC_n(x),\forall x$, and $D_n$ is 1-1.  So $\{D_n\}$ are
choice functions and each $D_n$ is 1-1.

\medskip
Define now the {\it index cocycle} $\pi_N:F|X_N\rightarrow S_N \;(=$
the symmetric group of $N$) by the formula:
\[
\pi_N(x,y)(k)=n\Leftrightarrow [C_k(x)]_E=[C_n(y)]_E
\]
(see Feldman-Sutherland-Zimmer \cite{FSZ}). Finally, we can define
$\sigma_N:[F|X_N]\times X_N\rightarrow S_N$ by
\[
\sigma_N(S,x)=\pi_N(x, S(x)).
\]

Since $S\in [F|X_N]$ is not a function but an equivalence class of
functions identified $\mu$-a.e., $\sigma_N$ again is to be
understood as an equivalence class of functions $(\sigma_N)_S (x)
=\sigma_N(S,x)$ identified $\mu$-a.e.  We again have the cocycle
identity:  For each $S,T\in [F|X_N]$,
\[
\sigma_N(ST,x)=\sigma_N(S,T(x))\sigma_N(T,x),
\]
for almost all $x\in X_N$.

Consider now the Hilbert space
\[
\h=\bigoplus_NL^2(X_N\times N),
\]
where $X_N\times N$ has the $\sigma$-finite measure $(\mu |X_n)\times\mu_N$, with
$\mu_N=$ the counting measure on $N$, and the unitary representation $\tau$ of
$[F]$ on $\h$ given by
\[\tau (S) (\bigoplus_N f_N) = \bigoplus_N g_N,
\]
\[
g_N(x,n)=f_N(S^{-1}(x),\sigma_N(S^{-1},x)(n)),
\]
for $(x,n)\in X_N\times N,f_N\in L^2(X_N\times N)$.  Clearly each
$L^2(X_N\times N)$ is invariant.

Notice that the representation $\tau$ is independent of the choice
functions $\{C_n^{(N)} \}$, up to unitary equivalence.

Consider the unit vector
\[
\xi_0=\bigoplus_N\chi_{X_N\times\{0\}}
\]
in $\h$.  Then for $S\in [F]$,
\begin{align*}
\langle\tau (S)(\xi_0),\xi_0\rangle &=\sum_N\int_{X_N\times
N}\xi_0(S^{-1}(x),\sigma_N
(S^{-1},x)(n))\xi_0(x,n)d\mu(x)d\mu_{N}(n)\\
&=\sum_N\mu (\{x\in X_N:\sigma_N(S^{-1},x)(0)=0\})\\
&=\sum_N\mu (\{x\in X_N:S^{-1}(x)Ex\})\\
&=\sum_N\mu (\{x\in X_N:xES(x)\})\\
&=\mu (\{x:S(x)Ex\})\\
&=\varphi_E(S).
\end{align*}
Thus the representation $\tau$ restricted to the closed span of
$\{\tau (S)(\xi_0): S\in [F]\}$ is the GNS representation of $[F]$
associated with $\varphi_E$.

If now $\Gamma$ is a countable group acting in a Borel way on $(X,
\mu)$ so that $E^X_{\Gamma} = F$, then, denoting by $\gamma$ also
the map $x\mapsto \gamma\cdot x$, the cocycle $\sigma_N$ restricts
to a cocycle, also denoted by $\sigma_N$, from $\Gamma\times X_N$ to
$S_N$: $\sigma_N (\gamma, x) = \pi_N(x, \gamma\cdot x)$. Similarly,
the representation $\tau$ restricts to a representation, also
denoted by $\tau$, of $\Gamma$ on $\h$.

We now characterize the condition on $E\subseteq F$ under which the representation
$\tau$ has an invariant non-0 vector.  First we note the following:
\medskip

\noindent{\bf Proposition 2.2}. {\it A vector $\xi$ is invariant under the
$\Gamma$-representation iff $\xi$ is invariant under the $[F]$-representation.}
\medskip

{\bf Proof}.  Suppose $\xi$ is invariant under the
$\Gamma$-representation, i.e., for $\xi =\bigoplus_N\xi_N$,
\[
\xi_N(x,n)=\xi_N(\gamma^{-1}\cdot x,\sigma_N(\gamma^{-1},x)(n)),
\]
for all $x\in X_N$, for all $\gamma\in\Gamma$ (neglecting as usual
null sets). Let now $S\in [F]$. Then for each $x\in X_N$, there is
$\gamma =\gamma_x\in\Gamma$ with $S^{-1}(x)=\gamma^{-1} \cdot x$.
Thus
\[
\tau (S)(\xi_N)(x,n)=\xi_N(S^{-1}(x),\sigma_N(S^{-1},x)(n))
\]
and $\sigma_N(S^{-1},x)(n)=k$, where
\begin{align*}
[C_n(x)]_E&=[C_k(S^{-1}(x))]_E\\
&=[C_k(\gamma^{-1}\cdot x)]_E,
\end{align*}
so $\sigma_N(\gamma^{-1},x)(n)=k=\sigma_N(S^{-1},x)(n)$, therefore
\begin{align*}
\tau (S)(\xi_N)(x,n)&=\xi_N(\gamma^{-1}\cdot x,\sigma_N(\gamma^{-1},x)(n))\\
&=\xi_N(x,n),
\end{align*}
i.e., $\xi_N$ and thus $\xi$ is also invariant under $\tau
(S)$.\hfill$\dashv$
\medskip

We now have:
\medskip

\noindent{\bf Proposition 2.3}. {\it The representation $\tau$ has
an invariant non-0 vector iff there is a Borel set $A\subseteq X$ of
positive measure which is $E$-invariant and $1\leq m <\infty$ such
that $[F|A:E|A]=m$, i.e., on some $E$-invariant Borel set of
positive measure $A$, there are exactly $m$ $E$-classes contained in
each $F|A$-class.} In particular, if $E$ is ergodic, $[F:E] <
\infty$.
\medskip

{\bf Proof}.  If such an $A$ exists, we can clearly assume that $A\subseteq X_N$
for some $N$.  Then let $B\subseteq X_N\times N$ be defined by
\[
(x,n)\in B\Leftrightarrow [C^N_n(x)]_E\subseteq A.
\]
Clearly for each $x\in X_N, B_x=\{n:(x,n)\in B\}$ has cardinality
$\leq m$, so if $\xi =\chi_B, \xi\in L^2(X_N\times N)$ and obviously
$\xi\neq 0$. Now we claim that $\xi$ is $\Gamma$-invariant, i.e.,
for $\gamma\in\Gamma$,
\[
\xi (x,n)=\xi (\gamma^{-1}\cdot x,\sigma_N(\gamma^{-1},x)(n)).
\]
This is clear, since
$[C^N_n(x)]_E=[C^N_{\sigma_N(\gamma^{-1},x)(n)}(\gamma^{-1}\cdot
x)]$, by the definition of $\sigma_N$.

Conversely, let $\xi\in\h$ be non-0 and $\Gamma$-invariant.  Clearly we can assume
that $\xi\in L^2(X\times N)$ for some $N$.  Now
\[
0<\int_{X_N}\sum_{n\in N}|\xi (x,n)|^2d\mu <\infty ,
\]
so for almost all $x\in X_N, \sum_{n\in N}|\xi (x,n)|^2<\infty$. Let
then $N_x= \{n\in N:|\xi (x,n)|$ is maximum among all $|\xi
(x,i)|,i\in N\}$.  Let $a_x$ be this maximum.  Then $N_x$ is finite,
provided that $a_x>0$.  Since $\xi$ is $\Gamma$-invariant
\[
\xi (x,n)=\xi (\gamma^{-1}\cdot x,\sigma_N(\gamma^{-1}, x)(n)),
\]
so $n\in N_x\Leftrightarrow\sigma_N(\gamma^{-1},x)(n)\in
N_{\gamma^{-1}\cdot x}, {\rm card}(N_x)={\rm
card}(N_{\gamma^{-1}\cdot x})$, and $a_x=a_{\gamma^{-1}\cdot x}$,
thus $x\mapsto {\rm card} (N_x), x\mapsto a_x$ are $F$-invariant.
Also as $$\int_{X_N}\sum_{n\in N}|\xi (x,n)|^2d\mu
>0 ,$$
$\{x\in X_N :a_x>0\}$ has positive measure. So fix $m > 0$ and a set
$Y\subseteq X_N$ of positive measure, which is $F$-invariant, and
for $x\in Y$ we have $a_x>0$ and $m= {\rm card} (N_x)$. Let then
\[
A=\bigcup\{[C^N_n(x)]_E:x\in Y,n\in N_x\}.
\]
Then $A$ is $E$-invariant, has positive measure and $[F|A:E|A]=m$.
\hfill$\dashv$

\medskip
Consider now the closed convex hull $C$ of
$\{\gamma\cdot\xi_0:\gamma\in\Gamma\}$, where $S\cdot\xi_0=\tau
(S)(\xi_0)$.  Since $\varphi_E(\gamma )=\langle\gamma\cdot\xi_0,
\xi_0\rangle$, we see that if
$\inf_{\gamma\in\Gamma}\varphi_E(\gamma )=\varphi^0_E
>0$, then $\langle\gamma\cdot\xi_0,\xi_0\rangle\geq\varphi^0_E,\forall\gamma\in
\Gamma$, so $\langle\eta ,\xi_0\rangle\geq\varphi^0_E,\forall\eta\in C$.  If then
$\xi$ is the unique element of least norm in $C$, we have $\langle\xi ,\xi_0\rangle
\geq\varphi^0_E$, and thus $\xi\neq 0$.  Clearly $\xi$ is invariant (under $\Gamma$
and thus $[F]$).  Thus we have
\medskip

\noindent{\bf Proposition 2.4}. {\it If
$\inf_{\gamma\in\Gamma}\varphi_E(\gamma )= \varphi^0_E>0$, then
there is a non-0 invariant vector for $\tau$.}
\medskip

{\bf (B)} We can conclude from 2.3 and 2.4 that if $\varphi^0_E>0$,
then there is an $E$-invariant set of positive measure $A$ such that
$[F|A:E|A]=m<\infty$. We can in fact obtain an estimate for such $m$
(and prove a somewhat stronger version).

\medskip
\noindent{\bf Theorem 2.5}. {\it Let $\Gamma$ be a countable group
and consider a measure preserving action of $\Gamma$ on $(X,\mu )$
with associated equivalence relation $F=E^X_\Gamma$.  Let
$E\subseteq F$ be a subequivalence relation.  Let $S,S'\in [F]$ and
assume that $\inf_{\gamma\in\Gamma}\varphi_E(S\gamma S')=c>0$. Then
there is an $E$-invariant Borel set $A$ of positive measure such
that $[F|A: E|A]=m\leq\tfrac{1}{c}$.  In particular, if
$c>\tfrac{1}{2}, F|A=E|A$.}

\medskip
{\bf Proof}.  Since $S\gamma S'=(SS')((S')^{-1}\gamma S')$, by
replacing the action $\gamma\cdot x$ of $\Gamma$ by the conjugate
action $\gamma *x = (S')^{-1}(\gamma\cdot S'(x))$, which also
induces $F$, we can assume that $S'={\rm id}$.  Thus we have
$\inf_{\gamma\in\Gamma}\varphi_E(S\gamma )=c>0$.  So $\langle
S\gamma\cdot\xi_0, \xi_0\rangle\geq c$, or
$\langle\gamma\cdot\xi_0,S^{-1}\cdot\xi_0\rangle\geq c,
\forall\gamma\in\Gamma$, thus if $C$ is the closed convex hull of
$\{\gamma\cdot \xi_0:\gamma\in \Gamma\}$, and $\xi$ the element of
least norm in $C$, $\xi$ is invariant for $\tau$ and $\langle\xi
,S^{-1}\cdot\xi_0\rangle =\langle S\cdot \xi ,\xi_0\rangle
=\langle\xi ,\xi_0\rangle\geq c$.

Now fix $\epsilon >0$ and let $\alpha_1,\dots ,\alpha_k\in[0,1],$
with $\sum^k_{i=1} \alpha_i=1$, and $\gamma_1,\dots
,\gamma_k\in\Gamma$ be such that if $\xi'=\sum^k_{i=1}
\alpha_i(\gamma_i\cdot\xi_0)$, then
$\Vert\xi'-\xi\Vert\leq\epsilon$.  Then, as $\xi$ is invariant, for
any $T\in [F]$ we have
\[
\langle T\cdot\xi',\xi_0\rangle =\langle T\cdot
(\xi'-\xi),\xi_0\rangle +\langle \xi ,\xi_0\rangle \geq c-\epsilon
\]
(note that $\langle T\cdot\xi',\xi_0\rangle$ is real). Thus for any
$T\in [F]$,
\[
\sum^k_{i=1}\alpha_i\langle T\gamma_i\cdot\xi_0,\xi_0\rangle\geq c-\epsilon
\]
or
\begin{equation}
\sum^k_{i=1}\alpha_i\varphi_E(T\gamma_i)\geq c-\epsilon .
\end{equation}
We will now use the following lemma.
\medskip

\noindent{\bf Lemma 2.6}. {\it Assume $E\subseteq F$ are countable,
measure preserving equivalence relations on $(X,\mu )$ and let
$n\geq 1$.  Then either there is an $E$-invariant Borel set
$A\subseteq X$ of positive measure such that every $F|A$-class
contains at most $n\ E|A$-classes or there are $T_0,\dots ,T_n \in
[F]$ such that $T_0 (x)=x$ and $[T_i(x)]_E\neq [T_j(x)]_E$, if
$i\not= j$.}
\medskip

Assuming the lemma, take $n=[\tfrac{1}{c}]$.  If the first case of
2.6 occurs, then the conclusion of the theorem immediately follows,
so it is enough to show that no such $T_0,\dots ,T_n$ exist.
Otherwise, apply (1) to $T_0,\dots ,T_n$ to get
\[
\sum^n_{j=0}\sum^k_{i=1}\alpha_i\varphi_E(T_j\gamma_i)\geq(n+1)(c-\epsilon
).
\]
But notice that
\[
\sum^n_{j=0}\varphi_E(T_jT)\leq 1,\forall T\in [F].
\]
This holds since
\[
\sum^n_{j=0}\varphi_E(T_jT)=\sum^n_{j=0}\mu (\{x:T_jT(x)Ex\})
\]
and the sets $\{x:T_jT(x)Ex\}, j=0,\dots ,n$ are pairwise disjoint. Thus
\[
(n+1)(c-\epsilon )\leq\sum^k_{i=1}\alpha_i=1,
\]
so, as $\epsilon$ is arbitrary, $n+1\leq \tfrac{1}{c}$, a
contradiction.

So it only remains to give the proof of 2.6.
\medskip

{\bf Proof of Lemma 2.6}. Assume first that $F$ is ergodic.

Consider then the ergodic decomposition of $E$. This is given by a
Borel map $\Sigma :X\rightarrow\e$, where $\e$ is the standard Borel
space of invariant, ergodic probability measures for $E$, such that
(i) $\Sigma$ is $E$-invariant; (ii) If $e\in\e$ and
$X_e=\Sigma^{-1}(\{e\})$, then $e(X_e)=1$ and $e$ is the unique
$E$-invariant probability measure on $X_e$; (iii) If $\Sigma_*\mu
=\mu_*$, then $\mu =\int ed\mu_*(e)$, i.e., $\mu (B)=\int
e(B)d\mu_*(e)$, for all Borel sets $B\subseteq X$.

Let $\e_0=$ atomic part of $\mu_*$, and put $\e_1=\e\setminus\e_0$.
Split $X_1= \bigcup_{e\in\e_1}X_e$ into $E$-invariant Borel sets
$X_1=A_0\sqcup\dots\sqcup A_n$, where $\mu (A_i)=\mu (A_j),\forall
i,j$,and let $\varphi_{i,j}\in [F|X_1]$ be such that
$\varphi_{i,j}(A_i)=A_j, 0\leq i,j\leq n$.  Now let $\psi_0,\dots,
\psi_n:\{0,\dots ,n\}\rightarrow \{0,\dots ,n\}$ be the bijections
defined by
\[
\psi_i(m)=(m+i)\mod (n+1).
\]
Then define $\varphi^{(1)}_i\in [F|X_1]$ by
\[
\varphi^{(1)}_i|A_m=\varphi_{m,\psi_i(m)}
\]
(so that $\varphi^{(1)}_i(A_m)=A_{\psi_i(m)}$).  Note that $\neg
\;\varphi^{(1)}_i(x) E \varphi^{(1)}_j(x)$, if $i\neq j$. Thus we
have found $\varphi^{(1)}_0,\dots , \varphi _n^{(1)}\in [F|X_1]$
with $\varphi^{(1)}_0(x)=x, \neg\; \varphi^{(1)}_i (x) E
\varphi^{(1)}_j (x)$, if $i\neq j$.

Consider now $e\in\e_0$, so that $\mu (X_e)>0$.  If
$[F|X_e:E|X_e]\leq n$, then $A=X_e$ satisfies the first alternative
of the lemma.  So we can assume that $[F|X_e:E|X_e] \geq n+1,\forall
e\in\e_0$. Since $E|X_e$ is ergodic, we can find $\varphi^e_0,\dots
,\varphi^e_n\in [F|X_e]$ with $\varphi^e_0(x)=x$ and $\neg\;
\varphi^e_i(x) E \varphi^e_j(x)$, if $i\neq j$ (see 2.1). Let
$\varphi^{(0)}_i=\bigcup_{e\in\e_0} \varphi^e_i$. Thus
$\varphi^{(0)}_i\in [F|X_0]$, where $X_0=\bigcup_{e\in\e_0}X_e
=X\setminus X_1$, and $\varphi^{(0)}_0(x)=x, \neg\; \varphi^{(0)}_i
(x) E  \varphi^{(0)}_j(x)$, if $i\neq j$.  Finally let
$T_i=\varphi^{(0)}_i\cup \varphi^{(1)}_i$.  This clearly works.

If $F$ is not ergodic, consider its ergodic decomposition and apply
the preceding argument to each piece of the ergodic decomposition.
\hfill$\dashv$

\medskip
We also have the following result concerning the ``proximity'' of
$E$ to $F$.
\medskip

\noindent{\bf Theorem 2.7}.  {\it Let $\Gamma$ be a countable group and consider
a measure preserving action of $\Gamma$ on $(X,\mu )$ with associated equivalence
relation $F=E^X_\Gamma$.  Let $E$ be a countable measure preserving equivalence
relation on $(X,\mu )$.  If $\epsilon >0$ is such that $\varphi_E(\gamma)\geq
1-\epsilon ,\forall\gamma\in\Gamma$, then $\varphi_E(S)\geq 1-4\epsilon ,\forall
S\in [F]$.}
\medskip

{\bf Proof}.  Since for $S\in[F]$, $\varphi_E(S)=\varphi_{E\cap
F}(S)$, we can assume that $E \subseteq F$.

In the earlier notation of Section 2, {\bf (A)} concerning the
representation $\tau$ of $[F]$, we have that
$\langle\gamma\cdot\xi_0,\xi_0\rangle\geq 1-\epsilon
,\forall\gamma\in\Gamma$ (where we put as before $\tau (S)(\xi
)=S\cdot\xi$).  If $\xi$ is the element of least norm in the closed
convex hull of $\{\gamma\cdot\xi_0 :\gamma\in\Gamma\}$, then $\xi$
is invariant for $\tau$ and $\langle\xi ,\xi_0\rangle \geq
1-\epsilon$, thus $\Vert\xi
-\xi_0\Vert^2\leq2(1-\langle\xi,\xi_0\rangle )\leq 2\epsilon$.  Thus
for any $S\in [F], \Vert \xi -S\cdot\xi_0\Vert^2 = \Vert S\cdot \xi
-S\cdot\xi_0\Vert^2\leq 2\epsilon$, so $\Vert
S\cdot\xi_0-\xi_0\Vert\leq 2\sqrt{2\epsilon}$, therefore
$2(1-\varphi_E(S))\leq 8\epsilon$, and so $\varphi_E (S)\geq
1-4\epsilon$.

\hfill$\dashv$

\medskip
{\bf (C)} We will next derive some consequences of the preceding
results. We refer to Gaboriau \cite{G2} for information about the
concept of {\it measure equivalence} (ME) introduced by Gromov.

\medskip
\noindent{\bf Lemma 2.8}.  {\it Let $\Gamma ,\Delta$ be two
countable groups and consider free, measuring preserving actions of
$\Gamma ,\Delta$ on $(X,\mu )$ with associated equivalence relations
$E=E^X_\Delta\subseteq F=E^X_\Gamma$.  If there is an $E$-invariant
Borel set $A\subseteq X$ such that every $F|A$-class contains at
most finitely many $E|A$-classes, then $\Gamma ,\Delta$ are {\rm
ME}.}

\medskip
{\bf Proof}.  By shrinking $A$ we can assume that there is $n$ such
that $[F|A:E|A]=n$. Fix Borel $T_1,\cdots ,T_n:A\rightarrow A$ such
that the $F|A$-class of $x$ is the union of the $E|A$-classes of
$T_i(x), i=1,\dots ,n$.

Let $\Omega =\{(x,y):x\in A,y\in X, xFy\}\subseteq F$ and consider
on $\Omega$ the $\sigma$-finite Borel measure $\nu |\Omega$, where
$\nu$ is the $\sigma$-finite Borel measure on $F$ given by $\nu
(B)=\int_X{\rm card}(B\cap F^x) d\mu (x)$, for every Borel
$B\subseteq F$.  Then $\Delta$ acts on $\Omega$ by: $\delta \cdot
(x,y)=(\delta\cdot x,y)$, since $A$ is $\Delta$-invariant, and
$\Gamma$ acts on $\Omega$ by: $\gamma\cdot (x,y)=(x,\gamma\cdot y)$.
These actions clearly preserve $\nu |\Omega$ and commute.  So it is
enough to show that each of the $\Delta ,\Gamma$ actions admits a
transversal (fundamental domain) of finite measure.  Let
$T:[A]_F\rightarrow A$ be a Borel map such that $T(y)Fy$. Then
$\Omega_1=\{(T_iT(y),y):y\in [A]_F, i=1,\dots ,n\}$ is a finite
measure transversal for the $\Delta$-action and
$\Omega_2=\{(x,x):x\in A\}$ is a finite measure transversal for the
$\Gamma$-action.\hfill$\dashv$

\medskip
\noindent{\bf Corollary 2.9}. {\it In the context of 2.8, if there
are $S,S'\in [F]$ such that $\inf_{\gamma\in\Gamma}\varphi_E(S\gamma
S')>0$, then $\Gamma ,\Delta$ are {\rm ME}.}

\medskip
\noindent{\bf Corollary 2.10}.  {\it Let $\Gamma$ be a countable
group and consider a free measure preserving action of $\Gamma$ on
$(X,\mu )$ with associated equivalence relation $F=E^X_\Gamma$.  If
$E\subseteq F$ is ergodic, treeable and there are $S,S'\in [F]$ such
that $\inf_{\gamma\in\Gamma}\varphi_E(S\gamma S')>0$, then $\Gamma$
has the Haagerup Approximation Property {\rm (HAP)}.}
\medskip

{\bf Proof}.  By Hjorth \cite{H3}, $E$ is given by a free action of
a group $\Delta$. Then $\Delta$ has the HAP (see, e.g., Gaboriau
\cite{G2}), and by 2.8 $\Gamma ,\Delta$ are ME, so $\Gamma$ has the
HAP (see, again Gaboriau \cite{G2}).\hfill$\dashv$
\medskip

Below, for an equivalence relation $F$, we denote by $[[F]]$ the set
of measure preserving bijections $\theta :{\rm dom} (\theta
)\rightarrow{\rm rng}(\theta)$ with ${\rm dom}(\theta ),{\rm
rng}(\theta )$ Borel sets and $\theta (x)F(x)$, for almost all
$x\in{\rm dom}(\theta )$.
\medskip

\noindent{\bf Lemma 2.11}. {\it Let $F$ be a countable, measure
preserving, ergodic equivalence relation on $(X,\mu )$ and let
$E\subseteq F$.  Let $X_\infty =\{x\in X: [x]_E$ is infinite\}.
Assume that there is an $E$-invariant Borel set $A$ of positive
measure such that $F|A=E|A$ (and thus $A\subseteq X_\infty$).  Then
for every $\epsilon >0$, there is $\theta \in [[F]]$ with $\mu ({\rm
dom}(\theta ))<\epsilon$ such that if $E\vee\theta$ is the
subequivalence relation generated by $E$ and $\theta$, then
$(E\vee\theta )|X_\infty =F|X_\infty$, a.e.  So if in addition $E$
is aperiodic (i.e., $X_\infty =X$), then $E\vee\theta =F$, a.e.

Moreover, if $F$ is given by an action of a finitely generated group $\Gamma$ and
$\{\gamma_1,\gamma_2,\dots ,\gamma_n\}$ is a set of generators for $\Gamma$, then
there are Borel sets $B_1,\dots ,B_n$ with $\mu (B_i)<\epsilon ,\forall i\leq n$,
such that if $E'=E\vee\gamma_1|B_1\vee\dots\vee\gamma_n|B_n$, then $E'|X_\infty
=F|X_\infty$.}
\medskip

{\bf Proof}.  Consider the ergodic decomposition of $E$, as in the proof of 2.6,
whose notation we keep below.

Since $\mu (A)>0$ and $E|A$ is ergodic, it follows that the measure
$\mu_*$ has atoms, and $A=X_{e_0}$, for some atom $e_0$ of $\mu_*$.
Let $\e_0$ be the set of atoms of $\mu_*$,
$\e_1=\{e\in\e\setminus\e_0:e$ is non-atomic\} and $\e_2=\e
\setminus (\e_0\cup\e_1)=\{\sigma\in\e :e$ is atomic\}.  Note that
$X_{\e_2}= \bigcup_{e\in\e_2}X_e=X_{\rm fin}=X\setminus X_\infty$.
We can clearly assume that $\mu_*(\e\setminus (\e_2\cup\{e_0\}))>0$,
otherwise there is nothing to prove.

Fix now $\mu (X_{e_0}) > \epsilon >0$.  If $e\not\in\e_2$, then $e$
is not atomic, so we can find $Y_e\subseteq X_e$ a Borel set with
$e(Y_e)=\epsilon /2$. Then, by ergodicity, $Y_e$ meets every
$E|X_e$-class. Let $Y=\bigcup_{e\not\in \e_2,e\neq e_0}Y_e$, so that
$0<\mu (Y)<\epsilon$. Then there is $\theta \in [[F]]$, with ${\rm
dom}(\theta )=Y, {\rm rng}(\theta )\subseteq X_{e_0}$.  We claim
that if $\bar E=E\vee\theta , \bar E|X_\infty =F|X_\infty$. To see
this note that if $y \in X_e$ for some $e\not\in \e_2, e\not= e_0$,
then there is $z\in Y_e$ with $yEz$. Thus $y\bar E \theta (z) \in
X_{e_0}$. So every $y \in X_{\infty}$ is $\bar E$-equivalent to an
element of $X_{e_0}$. Since $E|X_{e_0} = F|X_{e_0}$, we are done.

For the last assertion, decompose ${\rm dom}(\theta )=Y$ into
countably many Borel sets of positive measure
$\{Y_m\}^\infty_{m=1}$, so that there are words
$\{w_m\}^\infty_{m=1}$ in $\{\gamma_1,\dots ,\gamma_n\}$ with
$\theta |Y_m=w_m|Y_m$.  Say $w_m$ has length $k_m$.  Then find a
Borel set $Z_m\subseteq Y_m$ such that $\mu (Z_m)<\tfrac{1} {k_m}\mu
(Y_m)$ and for every $e\not\in\e_2,e\neq e_0$, if $e(Y_m\cap
X_e)>0$, then $e(Z_m\cap X_e)>0$ (again ignoring null sets for
$\mu_*$).  Let $Z=\bigcup_m Z_m$, so that $e(Z\cap X_e)>0$, for each
$e\neq\e_2, e\neq e_0$.  Let $\theta' =\theta |Z$.  Then as before
$(E\vee\theta')|X_\infty =F|X_\infty$.  Note now that if
$w_m=\gamma^{\pm 1}_{i_1}\gamma^{\pm 1}_{i_2}\dots\gamma^{\pm
1}_{k_m}$, then $\theta |Z_m=\theta'|Z_m$ is a composition of
$\gamma^{\pm 1}_{i_{k_m}}| Z_m,\gamma^{\pm
1}_{i_{k_m-1}}|\gamma^{\pm 1}_{i_{k_m}}(Z_m),\dots $  Thus ${\rm
graph}(\theta
|Z_m)\subseteq\gamma_1|B^{(m)}_1\vee\dots\vee\gamma_n|B^{(m)}_n =$
the equivalence relation generated by $\gamma_i|B^{(m)}_i, 1\leq
i\leq n $, where $B^{(m)}_i$ are Borel sets with $\mu
(B^{(m)}_i)\leq k_m\cdot\tfrac{1} {k_m}\cdot\mu (Y_m)=\mu
(Y_m),\forall i\leq n$. Let $B_i=\bigcup_mB^{(m)}_i$. Then $\mu
(B_i)\leq\sum_m\mu (Y_m)<\epsilon$ and $E\vee\theta'\subseteq E'$,
so $E'|X_\infty =F|X_\infty$. \hfill$\dashv$

\medskip
{\bf Remark 2.12}.  In the notation of 2.11, it is clear that, under
the same hypothesis, if $E$ is also ergodic, then $E=F$.

\medskip
Below we denote by $C_\mu (E)$ the cost of the countable, measure
preserving equivalence relation $E$ on $(X,\mu )$ (see Gaboriau
\cite{G1} or Kechris-Miller \cite{KM}).

\medskip
\noindent{\bf Corollary 2.13}. {\it Let $\Gamma$ be a countable
group and consider a measure preserving, ergodic action of $\Gamma$
on $(X,\mu )$ with associated equivalence relation $F =
E^{\Gamma}_X$. Let $E\subseteq F$ be a subequivalence relation. If
$E$ is aperiodic and for some $S,S'\in [F],\inf_{\gamma
\in\Gamma}\varphi_E(S\gamma S')>\tfrac{1}{2}$, then $C_\mu (F)\leq
C_\mu (E)$.} \vskip 5pt {\bf Proof}.  By 2.11 and 2.5.\hfill$\dashv$
\medskip

The following gives a quantitative version of (part of) 2.5.
\medskip

\noindent{\bf Lemma 2.14}. {\it Let $F$ be a measure preserving,
ergodic equivalence relation on $(X,\mu )$ and $E\subseteq F$ a
subequivalence relation. Assume that $\forall S\in [F] (\varphi_E
(S)>0$).  Then there is an $E$-invariant Borel set $A$ of positive
measure such that $E|A=F|A$.  Moreover, for some $S\in
[F],\{x:S(x)Ex\}\subseteq A$, so that $\mu (A)\geq\varphi_E(S)$.}
\medskip

{\bf Proof}.  Consider the ergodic decomposition of $E$ as in the proof of 2.6, whose
notation we keep below.

\medskip
{\bf Claim}.  $\mu_*$ has atoms, i.e., $\e_0 \not=\emptyset$.
\medskip

{\bf Proof}.  If $\mu_*$ is non-atomic, fix $A_*\subseteq\e$ with $\mu_*(A_*)=
\tfrac{1}{2}$.  If $\Sigma^{-1}(A_*)=A$, then $\mu (A)=\tfrac{1}{2}$ and $A$ is
$E$-invariant.  Let $T\in [F]$ be such that $T(A)=\sim A$.  Then $\varphi_E(T)=0$,
a contradiction.

\medskip
If $e\in\e_0$, then $\mu (X_e)>0$.  If for all $e\in\e_0, E|X_e\neq
F|X_e$, then as $E|X_e$ is ergodic, we can find $\varphi_e\in
[F|X_e]$ such that $\neg\; \varphi_e (x) E x,\forall x\in X_e$ (see
Remark 2.1).  If also $\mu_*(\sim\e_0)>0$, then we can split
$\sim\e_0$ into two pairwise disjoint sets of equal measure and thus
split $\sim X_{\e_0}=\sim\bigcup_{e\in\e_0}X_e$, into two sets of
equal measure $X_1,X_2$ which are $E$-invariant.  Then let $T_1\in
[F]$ be such that $T_1(X_1)=X_2, T_1(X_2)=X_1$.  Then $\neg \;
T_1(x) E x$ for $x\not\in X_{\e_0}$. If $T\in [F]$ is defined by
$T=(\bigcup_{e\in\e_0}\varphi_e )\cup (T_1|\sim X_{\e_0})$, then
clearly $\neg\; T(x) E x,\forall x$, a contradiction.

Thus we see that there must be some atom $e$ of $\mu_*$ with
$E|X_{e}=F|X_{e}$.

Enumerate in a sequence (finite or infinite) $\{e_0,e_1,\dots\}$ all
elements $e$ of $\e_0$ such that $E|X_e=F|X_e$ in such a way that
$\mu (X_{e_i})\geq\mu (X_{e_{i+1}})$. Put $A = X_{e_0}$. Let
$Y=\bigcup_iX_{e_i}, Z=\sim Y$.  We have seen that there is $S_1\in
[F|Z]$ with $\neg\; S_1(z)  E  z,\forall z\in Z$. Let for $n\geq 0$,
$\theta_{n+1} \in [[F]]$ be such that ${\rm dom} (\theta_{n+1}) =
X_{e_{n+1}}, {\rm rng}(\theta_{n+1})\subseteq X_{e_n}$, and let
$\theta^* = \bigcup_{n\geq 0} \theta_n$, so that $\theta^* :
\bigcup_{i>0} X_{e_i}\rightarrow \bigcup_i X_{e_i}, \theta^* \in
[[F]]$ and $\neg\; \theta^* (x) E x$. Let $\theta^{**} \in [[F]]$ be
such that ${\rm dom}( \theta^{**})=X_{e_0}$ and ${\rm
rng}(\theta^{**})=\bigcup_i X_{e_i}\setminus
\theta^*(\bigcup_{i>0}X_{e_i})$.  Put $S_2=\theta^*\cup\theta^{**}$,
so that $S_2 \in [F|Y]$ and $\neg\; S_2(y) E y$, if $y\not\in
X_{e_0}$.  Then if $S=S_1\cup S_2, S\in [F]$ and
$\{x:S(x)Ex\}\subseteq X_{e_0}=A$.\hfill$\dashv$
\medskip

\noindent{\bf Corollary 2.15}. {\it Let $\Gamma$ be a countable
group and consider a measure preserving, ergodic action of $\Gamma$
on $(X,\mu )$ with associated equivalence relation $F=E^X_\Gamma$.
Let $E\subseteq F$ be a subequivalence relation. If
$\varphi^0_E=\inf_{\gamma\in\Gamma}\varphi_E(\gamma) >\tfrac{3}{4}$,
then there is an $E$-invariant Borel set $A$ with $E|A=F|A$ such
that $\mu (A)\geq 4\varphi^0_E-3$.}

\medskip
{\bf Proof}.  By 2.14 and 2.7.\hfill$\dashv$

\section{Epstein's co-inducing construction}

{\bf (A)} We will next study some properties of a co-inducing construction of
Epstein \cite{E}.

We first describe this construction.  Fix a standard measure space
$(X,\mu )$, a countable, measure preserving equivalence relation $F$
on $(X,\mu )$ and a subequivalence relation $E\subseteq F$ such that
there is a fixed number $N\in\{1,2,3,\dots ,\aleph_0\}$ of
$E$-classes in each $F$-class.  This is the case, for example, if
$F$ is ergodic.  Fix choice functions $\{C_n\}_{n\in N}$, where we
identify $N$ here with $\{0,1,\dots ,N-1\}$, if $N$ is finite, and
with $\bbN$ if it is infinite, as in Section 2, {\bf (A)}, and let
$\pi :F\rightarrow S_N$ (= the symmetric group of $N$) be the index
cocycle given by the formula:
\[
\pi (x , y)(k)=n\Leftrightarrow [C_k(x)]_E=[C_n(y)]_E.
\]

Now assume that $E$ as above is induced by a {\it free} action $a_0$
of a countable group $\Delta$. Then we can define $\bar\delta
:F\rightarrow\Delta^N$ by
\[
\bar\delta (x,y)_n\cdot C_{\pi (x,y)^{-1}(n)}(x)=C_n(y).
\]
The group $S_N$ of permutations of $N$
acts on $\Delta^N$ by shift $(\pi\cdot\bar\delta )_n=
\bar\delta_{\pi^{-1}(n)}$, so we can consider the semi-direct product
$S_N\ltimes\Delta^N$, whose multiplication is defined by:
\[
(\pi_1,\bar\delta_1)(\pi_2,\bar\delta_2)=(\pi_1\pi_2,\bar\delta_1(\pi_1\cdot\bar\delta_2)).
\]
It is easy to check that
\[
\rho (x,y)=(\pi (x,y),\bar\delta (x,y))
\]
is a Borel cocycle
\[
\rho :F\rightarrow S_N\ltimes\Delta^N.
\]

Now given any measure preserving action $a$ of $\Delta$ on a standard
measure space $(Y,\nu )$, we can define a measure preserving action of
$S_N\ltimes\Delta^N$ on $(Y^N,\nu^N)$ by
\[
((\pi ,\bar\delta )\cdot\bar y)_n=\bar\delta_n\cdot\bar y_{\pi^{-1}(n)}.
\]
Then we can define a {\it near-action} of $[F]$ on $(X\times Y^N,
\mu \times \nu^N)$, i.e., a continuous homomorphism of $[F]$, with
the uniform topology, into the automorphism group ${\rm Aut}(X\times
Y^N, \mu \times \nu^N)$, with the weak topology, by letting $S\in
[F]$ act on $X\times Y^N$ as a skew product via $\rho$, namely

\begin{align*}
S\cdot (x,\bar y) &=(S(x),\rho (x, S(x))\cdot\bar y)\\
&=(S(x),(n\mapsto\bar\delta (x,S(x))_n\cdot\bar y_{\pi (x,
S(x))^{-1}(n)})).
\end{align*}
(For information about near-actions, see Glasner-Tsirelson-Weiss
\cite{GTW}, and concerning the uniform and weak topologies, see
Kechris \cite{K}.)

In particular, if $c_0$ is a measure preserving action of a
countable group $\Lambda$ on $(X,\mu)$ with $E_{\Lambda}^X\subseteq
F$, then $\lambda \in \Lambda$ gives rise to an element $x\mapsto
\lambda\cdot x$ of $[F]$, so we can define $$ \pi(\lambda, x) = \pi
(x, \lambda\cdot x), \bar \delta (\lambda, x) =
\bar\delta(x,\lambda\cdot x),$$
$$\rho (\lambda, x) = (\pi (\lambda, x),\bar\delta (\lambda, x))$$
and then $\rho: \Lambda\times X\rightarrow S_N\ltimes \Delta^N$ is a
Borel cocycle. The restriction of the near action of $[F]$ on
$(X\times Y^N, \mu \times \nu^N)$ gives then a measure preserving
action of $\Lambda$ on $(X\times Y^N, \mu \times \nu^N)$, which is
the skew product
$$c=c_0\ltimes_\rho (Y^N, \mu^N)$$defined by

\begin{align*}
\lambda\cdot (x,\bar y) &=(\lambda\cdot x,\rho (\lambda ,x)\cdot\bar y)\\
&=(\lambda\cdot x,(n\mapsto\bar\delta (\lambda ,x)_n\cdot\bar y_{\pi
(\lambda , x)^{-1}(n)}).
\end{align*}

Fix now a Borel action $b_0$ of a countable group $\Gamma$ on
$(X,\mu)$ with $F=E_{\Gamma}^X$. Then applying the above to $\Lambda
= \Gamma, c_0 = b_0$, we associate to each measure preserving action
$a$ of $\Delta$ on $(Y,\nu)$ a measure preserving action $b$ of
$\Gamma$ on $(X\times Y^N,\mu\times\nu^N)$, relative to the fixed
pair $(a_0,b_0)$ of the actions of $\Gamma ,\Delta$, resp., on $X$
(and the choice of $\{C_n\}$ -- but it is not hard to check that
this action is independent of the choice of $\{ C_n \}$, up to
isomorphism). We call this the {\it co-induced action of $a$, modulo
$(a_0,b_0)$}, in symbols:
\[
b={\rm C}{\rm Ind}(a_0,b_0)^\Gamma_\Delta (a).
\]
We can view this as an operation from the space $A(\Delta ,Y,\nu )$ of
measure preserving actions of $\Delta$ on $(Y,\nu )$ (see, e.g., Kechris
\cite{K}) to the space $A(\Gamma ,X\times Y^N,\mu\times\nu^N)$.

By applying the preceding to $\Lambda = \Delta, c_0 = a_0$, we also
have a measure preserving action $a'$ of $\Delta$ on $(X\times
Y^N,\mu\times\nu^N)$. Clearly this action gives a subequivalence
relation of the equivalence relation given by $b$. We note that
$b_0$ is a factor of $b$ via $(x,\bar y)\mapsto x$, $a$ is a factor
of $a'$ via $(x,\bar y)\mapsto \bar y_0$ (recall here that
$C_0(x)=x$, so that $\pi (\gamma ,x)(0)=0$ iff $\gamma\cdot x E x$)
and finally $a_0$ is a factor of $a'$ via $(x,\bar y)\mapsto x$. In
particular, if $b_0$ is free, so is $b$, and $a'$ is always free.

Finally, for further reference, we note that the map $a\mapsto {\rm C Ind}
(a_0,b_0)^\Gamma_\Delta (a)$ is a continuous map from $A(\Delta ,Y,\nu )$
to $A(\Gamma ,X\times Y^N,\mu\times\nu^N)$, where each is equipped with
the weak topology (see Kechris \cite{K} for its definition).
\medskip

{\bf (B)} We will now study some connections between ergodicity
properties of an action and its co-induced action.  {\it In the
notation of {\bf (A)}, if $a_0$ is ergodic we can choose the choice
functions to be 1-1 and it will be assumed in this case that the
co-inducing construction is done with such choice functions.}

\medskip
\noindent{\bf Proposition 3.1}. {\it In the notation of {\rm{\bf
(A)}} above, if $b_0$ is free, mixing and $a_0$ is ergodic, then
\[
a\text{ is mixing }\Rightarrow b={\rm CInd}(a_0,b_0)^\Gamma_\Delta (a)
\text{ is mixing}.
\]}
\indent {\bf Proof}.  Assume that $b_0,a$ are mixing, $b_0$ is free,
and $a_0$ is ergodic and consider the action $b$ of $\Gamma$ on
$(X\times Y^N,\mu \times\nu^N)$.  Then $\Gamma$ acts on $L^2(X\times
Y^N,\mu\times\nu^N)$ by $\gamma\cdot f(x,\bar y)=f(\gamma^{-1}\cdot
(x,\bar y))$ and it is enough to show that for $f,g\in L^2(X\times
Y^N,\mu\times\nu^N),\int (\gamma^{-1}\cdot f)g\rightarrow (\int
f)(\int g)$ as $\gamma\rightarrow \infty$.  Without loss of
generality, we can assume that $f(x,\bar y)=f_0 (x)F_0(\bar
y_0)\dots F_m(\bar y_m),g(x,\bar y)=g_0(x)G_0(\bar y_0)\dots
G_m(\bar y_m)$, for bounded $f_0,g_0:X\rightarrow\bbC , F_i,G_i:Y
\rightarrow\bbC$.  Note that since $b_0$ is mixing, $\int
f_0(\gamma\cdot x)g_0(x)d\mu (x)\rightarrow (\int f_0)(\int g_0)$.
Now we have
\[
\int (\gamma^{-1}\cdot f)g=\int f_0(\gamma\cdot x)g_0(x)[\int \prod^m_{i=0}
F_i(\bar\delta (\gamma ,x)_i\cdot\bar y_{\pi (\gamma ,x)^{-1}(i)})
G_i(\bar y_i)d\bar y]dx.
\]
Note that if $\{q_n\},q_n:Z\rightarrow\bbC ,\{r_n\},r_n:Z\rightarrow\bbC$,
are uniformly bounded, where $Z$ is a probability space, $q_n(z)\rightarrow
a,\forall z$, and $\int r_n\rightarrow b$, then
\[
\int q_n(z)r_n(z)=\int
(q_n(z)-a)r_n(z)+a\int r_n(z)\rightarrow ab,
\]
by Lebesgue Dominated Convergence.  So it is enough to show that
{\it for each fixed} $x\in X$,
\begin{gather}
\int \left [\prod^m_{i=0}F_i(\bar\delta (\gamma ,x)_i\cdot\bar
y_{\pi (\gamma , x)^{-1}(i)})G_i (\bar y_i)\right ]d\bar
y\rightarrow\prod^m_{i=0}\left (\int F_i \right )\left (\int
G_i\right ),\tag{$*$}
\end{gather}
as $\gamma\rightarrow\infty$.

Fix then $x\in X$ and put
\[
S_\gamma =\{(i,\pi (\gamma ,x)^{-1}(i)): i\leq m,\pi (\gamma
,x)^{-1}(i)\leq m\}.
\]
For each $S\subseteq \{0,\dots ,m\}^2$, let
\[
\Gamma_S=\{\gamma\in\Gamma :S_\gamma =S\}.
\]
Then $\Gamma =\bigsqcup_S\Gamma_S$ is a finite partition of $\Gamma$, so it
is enough to show that for each fixed $S$, with $\Gamma_S$ infinite, $(*)$
holds as $\gamma\rightarrow\infty ,\gamma\in\Gamma_S$.

For such $S$ and $\gamma\in\Gamma_S$, let
\[
I_\gamma =\{i\leq m:\pi (\gamma ,x)^{-1}(i)\leq m\}
\]
and let
\[
\rho_\gamma (i)=\pi (\gamma ,x)^{-1}(i),
\]
for $i\in I_\gamma$.  Thus graph$(\rho_\gamma )=S$.

Then
\begin{align*}
\int \biggl [ (\prod^m_{i=0}& F_i(\bar\delta (\gamma ,x)_i\cdot\bar y_{\pi (\gamma
,x)^{-1}(i)})G_i(\bar y_i)\biggr ]d\bar y\\
&=\int \left [\prod_{i\in I_\gamma}F_i(\bar\delta (\gamma ,x)_i\cdot\bar
y_{\rho_\gamma (i)})G_{\rho_\gamma (i)}(\bar y_{\rho_\gamma (i)})\right ]\\
&\quad\quad \left [\prod_{i\not\in I_\gamma} F_i (\bar\delta (\gamma
,x)_i \cdot \bar y_{\pi (\gamma
,x)^{-1}(i)})\prod_{i\not\in\rho_\gamma (I_\gamma )} G_i(\bar
y_i)\right ]d\bar y.
\end{align*}
Noticing that if $i\not\in I_\gamma ,i\leq m$, then $\pi (\gamma
,x)^{-1} (i)>m$, so that $$\rho_\gamma (I_\gamma ), \{\pi (\gamma
,x)^{-1}(i): i\not\in I_\gamma\}, \{0,\dots ,m\}\setminus\rho_\gamma
(I_\gamma )$$are pairwise disjoint, and applying independence, we
see that the above integral is equal to
\[
\left [\prod_{i\not\in I_\gamma}\int F_i\right ]\left [
\prod_{i\not\in\rho_\gamma (I_\gamma )}
\int G_i\right ]\
\left [\int\prod_{i\in I_\gamma} F_i(\bar\delta (\gamma ,x)_i\cdot
\bar y_{\rho_\gamma (i)})G_{\rho_\gamma (i)}(\bar y_{\rho_\gamma (i))}
)d\bar y\right ].
\]
Thus for each $i\in I_\gamma$, if $j=\rho_\gamma (i)$, it is enough, by
independence again, to show that
\[
\lim_{\gamma\in\Gamma_S,\gamma\rightarrow\infty}\int F_i(\bar\delta (\gamma
,x)_i\cdot\bar y_j)G_j(\bar y_j)d\bar y=\left (\int F_i\right )\left (\int G_j\right )
\]
or equivalently
\[
\int F_i(\bar\delta (\gamma ,x)_i\cdot y)G_j(y)dy\rightarrow \left (\int F_i\right )
\left (\int G_j\right )
\]
as $\gamma\rightarrow\infty ,\gamma\in \Gamma_S$.  Using that the
action $a$ of $\Delta$ on $Y$ is mixing, it is then enough to show
that for each $i$,
\[
\gamma\rightarrow\infty\Rightarrow\bar\delta (\gamma ,x)_i\rightarrow
\infty .
\]
Otherwise, there is a finite $K\subseteq\Delta$ such that for
infinitely many $\gamma\in\Gamma ,\bar\delta (\gamma ,x)_i\in K$.
Now $\bar\delta (\gamma , x)_i\cdot C_j(x)=C_i(\gamma \cdot x)$, so
$C_i(\gamma \cdot x)$ takes only finitely many values for infinitely
many $\gamma\in\Gamma$, contradicting the fact that the
$\Gamma$-action on $X$ is free and $C_i$ is 1-1.\hfill$\dashv$

\medskip
There is another condition, concerning the ``smallness'' of $E$ in
$F$, that actually guarantees that the co-induced action $b$ is
mixing for {\it any} $a$ (mixing or not).

In the context of {\bf (A)}, let for each $\gamma\in\Gamma, k,n\in
N$,
\begin{align*}
A^{k,n}_E(\gamma )&=\{x:\pi (\gamma ,x)(k)=n\}.\\
&=\{x:C_k(x)EC_n(\gamma\cdot x)\}.
\end{align*}
Thus $A^{0,0}_E(\gamma )=\{x:\gamma\cdot x Ex\}$.  Put
\[
\varphi^{k,n}_E(\gamma )=\mu (A^{k,n}_E(\gamma )).
\]
Clearly $\varphi^{0,0}_E=\varphi_E$.

\medskip
\noindent{\bf Lemma 3.2}. {\it If $\varphi_E(\gamma )\rightarrow 0$
as $\gamma \rightarrow\infty$, then for any $k,n,
\varphi^{k,n}_E(\gamma )\rightarrow 0$ as
$\gamma\rightarrow\infty$.}
\medskip

{\bf Proof}.  We have $\varphi^{k,n}_E(\gamma )=\mu
(\{x:C_k(x)EC_n(\gamma\cdot x)\})$.  Put $S=C_k,T=C_n$.  There is a
partition $X=\bigsqcup_{i\in\bbN} A_i$, and $\gamma_i\in\Gamma$ such
that $S=\bigsqcup_i\gamma_i|A_i$.  Similarly there is a partition
$X=\bigsqcup_{j\in\bbN}B_j$ and $\delta_j\in\Gamma$ such that
$T=\bigsqcup_j\delta_j|B_j$.

Assume $\varphi_E(\gamma )\rightarrow 0$ as
$\gamma\rightarrow\infty$, and let $\epsilon >0$.  We will find a
finite set $F\subseteq\Gamma$ such that for $\gamma\not\in
F,\varphi^{k,n}_E(\gamma )<\epsilon$.  First find $J_0$ such that
$\sum_{j\geq J_0}\mu (B_j)<\epsilon /3$.  Then fix $I_0$ such that
$\sum_{i\geq I_0}\mu (A_i)<\epsilon /3|J_0|$.  Since for each $i,j
,\gamma_i\gamma^{-1} \delta_j\rightarrow\infty$ as
$\gamma\rightarrow\infty$, there is a finite set $F\subseteq\Gamma$
such that for $\gamma\not\in F,\sum_{i<I_0,j<J_0}\varphi_E
(\gamma_i\gamma^{-1}\delta^{-1}_j)<\epsilon /3$.  We show that if
$\gamma\not\in F$, then $\varphi^{k,n}_E(\gamma )<\epsilon$.

We have for any $\gamma$:
\begin{align*}
\varphi^{k,n}_E(\gamma )&=\mu (\{x:S(x)ET(\gamma\cdot x)\})\\
&=\sum_{i,j}\mu (\{x:x\in A_i\wedge x\in\gamma^{-1}\cdot
B_j\wedge\gamma_i\cdot x
E\delta_j\gamma\cdot x\})\\
&=\sum_{i,j}\mu (\{x:\gamma^{-1}\delta^{-1}_j\cdot x\in A_i\wedge\gamma^{-1}
\delta^{-1}_j\cdot x\in\gamma^{-1}\cdot B_j\\
&\quad\quad\quad\wedge\gamma_i\gamma^{-1}\delta^{-1}_j\cdot xEx\})\\
&=\sum_{i,j}\mu (\{x:\gamma^{-1}\delta^{-1}_j\cdot x\in A_i\wedge\delta^{-1}_j\cdot
x\in B_j\\
&\quad\quad\quad\wedge\gamma_i\gamma^{-1}\delta^{-1}_j\cdot xEx\}).
\end{align*}
Let
\[
A^{(\gamma )}_{i,j}=\{x:\gamma^{-1}\delta^{-1}_j\cdot x\in A_i\wedge\delta^{-1}_j
\cdot x\in B_j\}.
\]
Then
\begin{align*}
\mu (A^{(\gamma )}_{i,j})&=\mu (\{x:\gamma^{-1}\cdot x\in A_i\wedge x\in B_j\})\\
&=\mu (\{x:x\in A_i\wedge\gamma\cdot x\in B_j\}),
\end{align*}
so $\sum_{j\geq J_0}\sum_i\mu (A^{(\gamma )}_{i,j})=\sum_{j\geq
J_0}\mu (B_j)<\epsilon /3$. Also $\sum_{j<J_0}\sum_{i\geq I_0}\mu
(A^{(\gamma )}_{i,j})\leq\break \sum_{j<J_0}\sum_{i\geq I_0}\mu
(A_i)<\epsilon /3$.  So it follows that
\[
\varphi^{k,n}_E(\gamma )\leq [\sum_{i<I_0}\sum_{j<J_0}\varphi_E(\gamma_i\gamma^{-1}
\delta^{-1}_j)]+2\epsilon /3,
\]
thus if $\gamma\not\in F$,
\[
\varphi^{k,n}_E(\gamma )<\epsilon .
\]
\hfill$\dashv$

\medskip
We denote below by $i$ the trivial action of $\Delta$ on $(Y,\nu
):\delta \cdot y=y$.  We now have:

\medskip
\noindent{\bf Theorem 3.3}.  {\it In the notation of {\rm {\bf(A)}}
above, and assuming that $b_0$ is mixing, the following are
equivalent:

(i) $\varphi_E(\gamma )\rightarrow 0$ as $\gamma\rightarrow \infty
$,

(ii) ${\rm CInd}(a_0,b_0)^\Gamma_\Delta (i)$ is mixing,

(iii) $\exists a\in A(\Delta ,Y,\nu )$ {\rm (}$a$ is not ergodic and
${\rm CInd} (a_0,b_0)^\Gamma_\Delta (a)$ is mixing{\rm )},

(iv) $\forall a\in A(\Delta ,Y,\nu )({\rm
CInd}(a_0,b_0)^\Gamma_\Delta (a)$ is mixing{\rm )}.}
\medskip

{\bf Proof}.   Clearly (ii) $\Rightarrow$ (iii) and (iv) $\Rightarrow$ (ii).

(i) $\Rightarrow$ (iv): By 3.2 we have that $\varphi^{k,n}_E(\gamma
)\rightarrow 0, \forall k,n$.  Then going over the proof of 3.1 and
keeping its notation, we see that for $x\not\in\bigcup_{1\leq
i,j\leq m}\{x:\pi (\gamma ,x)^{-1} (i)=j\}=\bigcup_{1\leq i,j\leq
m}A^{j,i}_E(\gamma )=A^{(m)}_E(\gamma )$,
\[
\int[\prod^m_{i=0}F_i(\bar\delta (\gamma ,x)_i\cdot\bar y_{\pi
(\gamma ,x)^{-1}(i)})G_i(\bar y_i)]d\bar y=\prod^m_{i=1}(\int
F_i)(\int G_i),
\]
by independence.  Thus, for some bounded $H(x)$,
\begin{align*}
\int (\gamma^{-1}\cdot f)g&=\int_{A^{(m)}_E(\gamma )}H(x)dx+\int_{\sim A^{(m)}_E(\gamma)}
f_0(\gamma\cdot x)g_0(x)[\prod^m_{i=0}(\int F_i)(\int G_i)]dx\\
&\rightarrow (\int f)(\int g)
\end{align*}
as $\gamma\rightarrow\infty$, since $\mu (A^{(m)}_E(\gamma ))\rightarrow 0$
and $b_0$ is mixing.

(iii) $\Rightarrow$ (i): Fix such an action $a$ and a set
$B\subseteq Y$ with $0<p=\mu (B)<1$, which is invariant under this
action. We will show that $\varphi_E(\gamma)\rightarrow 0$. Put
$B^{(0)} =\{(x,\bar y):\bar y_0\in B\}$. Since the co-induced action
is mixing, we have that $(\mu\times\nu^N)(\gamma\cdot B^{(0)}\cap
B^{(0)})\rightarrow (\mu\times\nu^N)(B^{(0)})\cdot (\mu\times\nu^N)
(B^{(0)})=p^2$.  Now $\gamma\cdot B^{(0)}=\{\gamma\cdot (x,\bar y):
(x,\bar y)\in B^{(0)}\}=\{\gamma\cdot (x,\bar y): \bar y_0\in B\}$,
so $\gamma\cdot B^{(0)}\cap B^{(0)}=\{\gamma\cdot (x,\bar y):\bar
y_0\in B\wedge (\rho (\gamma ,x)\cdot \bar y)_0\in B\}
=\{\gamma\cdot (x,\bar y):\bar y_0\in B\wedge\bar\delta (\gamma\cdot
x)_0\cdot\bar y_{\pi (\gamma ,x)^{-1}(0)}\in B\} =\{\gamma\cdot
(x,\bar y):\bar y_0\in B\wedge \bar y_{\pi (\gamma ,x)^{-1}(0)}\in
B\} =\{\gamma\cdot (x,\bar y):\pi (\gamma ,x)(0)=0 \wedge\bar y_0\in
B\}\cup \{\gamma\cdot (x,\bar y):\pi (\gamma ,x)(0)\neq 0 \wedge\bar
y_0\in B\wedge \bar y_{\pi (\gamma ,x)^{-1}(0)}\in B\}$. So, by
Fubini,
\begin{align*}
(\mu\times\nu^N)(\gamma\cdot B^{(0)}\cap B^{(0)})&=\mu
(A^{0,0}_E(\gamma
))\cdot\mu (B)+(1-\mu (A^{0,0}_E(\gamma ))\cdot\mu (B)^2\\
&=p\mu (A^{0,0}_E(\gamma ))+p^2(1-\mu (A^{0,0}_E(\gamma ))).
\end{align*}
Since $(\mu\times\nu^N)(\gamma\cdot B^{(0)}\cap B^{(0)})\rightarrow
p^2$ and $0<p<1, \mu (A^{0,0}_E(\gamma )) =
\varphi_E(\gamma)\rightarrow 0$.

\hfill$\dashv$

\medskip
It follows that if for some $k,n,\varphi^{k,n}_E(\gamma
)\not\rightarrow 0$, as $\gamma\rightarrow\infty$, then for every
$a\in A(\Delta ,Y,\nu )$, if ${\rm CInd}(a_0,b_0)^\Gamma_\Delta (a)$
is mixing, then $a$ is ergodic. By strengthening the hypothesis, we
can obtain the following stronger conclusion.
\medskip

\noindent{\bf Proposition 3.4}. {\it In the notation of {\rm {\bf
(A)}} above, if for some $k,n$, we have
$\inf_{\gamma\in\Gamma}\varphi^{k,n}_E(\gamma )>0$, then for any
$a\in A(\Delta ,Y,\nu )$, if $b={\rm CInd}(a_0,b_0)^\Gamma_\Delta
(a)$ is ergodic, then $a$ is ergodic.}
\medskip

{\bf Proof}.  Assume that $a$ is not ergodic, in order to show that $b$
is not ergodic.  Let $f\in L^2_0(Y)=\{f\in L^2(Y):\int f=0\}$ be real
such that $\Vert f\Vert_2=1$ and $\delta\cdot f=f,\forall\delta\in\Delta$.
Let $f^{(k)},f^{(n)}\in L^2_0(X\times Y^N)$ be defined by $f^{(k)}(x,
\bar y)=f(\bar y_k),f^{(n)}(x,\bar y)=f(\bar y_n)$.  Then for the
$\Gamma$-action on $L^2_0(X\times Y^N)$,
\begin{align*}
\langle\gamma^{-1}\cdot f^{(n)},f^{(k)}\rangle &=\iint f^{(n)}(\gamma\cdot
(x,\bar y))f^{(k)}(x,\bar y)dxd\bar y\\
&=\iint f^{(n)}(\gamma\cdot x,(n\mapsto\bar\delta (\gamma ,x)_n\cdot
\bar y_{\pi (\gamma ,x)^{-1}(n)}))f^{(k)}(x,\bar y)dxd\bar y\\
&=\iint f(\bar\delta (\gamma ,x)_n\cdot\bar y_{\pi (\gamma ,x)^{-1}(n)}
)f(\bar y_k)dxd\bar y\\
&=\int_{A^{k,n}_E(\gamma )}\left [\int f(\bar\delta (\gamma ,x)_n\cdot \bar y_k)
f(\bar y_k)d\bar y\right ]dx\\
&\quad\quad +\int_{\sim A^{k,n}_E(\gamma )}\left [\int f(\bar\delta (\gamma ,x)_n\cdot
\bar y_{\pi (\gamma ,x)^{-1}(n)})f(\bar y_k)d\bar y\right ]dx\\
&=\int_{A^{k,n}_E(\gamma )}\left (\int f^2\right )dx +\int_{\sim A^{k,n}_E(\gamma )}
\left (\int f\right )^2dx\\
&=\mu (A^{k,n}_E(\gamma ))=\varphi^{k,n}_E(\gamma ).
\end{align*}
Thus for some $c>0,\langle\gamma^{-1}\cdot
f^{(n)},f^{(k)}\rangle\geq c>0$.  If $K$ is the closed convex hull
of $\{\gamma\cdot f^{(n)}:\gamma \in\Gamma \}$, then $\langle\xi
,f^{(k)}\rangle\geq c,\forall\xi \in K$, so $0\not\in K$.  If
$\xi_1$ is the unique element of least norm in $K$, clearly
$0\neq\xi_1\in L^2_0(X\times Y^\bbN )$, and $\xi_1$ is
$\Gamma$-invariant, so the action $b$ is not ergodic. \hfill$\dashv$

\medskip
Let $b_0$ be a free action of a countable group $\Gamma$ on $(X,\mu
)$ and let $a_0$ be a free action of a countable group $\Delta$ on
$(X,\mu )$, so that $E=E^X_\Delta\subseteq F=E^X_\Gamma$.  If there
exist $S,S'\in [F]$ with $\inf_{\gamma\in\Gamma}\varphi_E(S\gamma
S')>0$, then, by 2.9, $\Gamma ,\Delta$ are ME.  In particular, if
$a_0$ is ergodic, so that we can take the choice functions $\{C_n\}$
to be in $[F]$, we note that $\varphi^{k,n}_E(\gamma )=\mu
(\{x:C_k(x)EC_n(\gamma\cdot x)\})=\mu (\{x:C_n\gamma
C^{-1}_k(x)Ex\})=\varphi_E (C_n\gamma C^{-1}_k)$, thus we have:
\medskip

\noindent{\bf Corollary 3.5}. {\it Let $b_0$ be a free measure
preserving action of $\Gamma$ on $(X,\mu )$ and let $a_0$ be a free,
ergodic action of $\Delta$ on $(X,\mu )$ with $E=E^X_\Delta\subseteq
F=E^X_\Gamma$.  If for some $k,n,
\inf_{\gamma\in\Gamma}\varphi^{k,n}_E(\gamma )>0$, then $\Gamma
,\Delta$ are {\rm ME}.  In particular, $\Gamma$ has property {\rm
(T)} {\rm (}resp., {\rm HAP}{\rm )} iff $\Delta$ has property {\rm
(T)} {\rm (}resp., {\rm HAP}{\rm )}.}

\medskip
{\bf Remark 3.6}. In the context of 3.5, when $b_0$ is also mixing
and $\Gamma$ does not have the HAP, one can show that $\Delta$ does
not have the HAP by the following alternative argument, which may be
of some independent interest.

We will use the following characterization of groups
with HAP, see Kechris \cite{K}, 12.7.
Below ERG$(\Gamma ,X,\mu )$ is the set of ergodic actions of $\Gamma$ on
$(X,\mu )$, and MIX$(\Gamma ,X,\mu )$ the set of mixing actions.  We
consider these as subspaces of the space of actions $A(\Gamma ,X,\mu )$
with the weak topology.

\medskip
\noindent{\bf Theorem 3.7}. {\it Let $\Gamma$ be an infinite
countable group.  Then the following are equivalent:

(i) $\Gamma$ does not have the {\rm HAP},

(ii) $\overline{{\rm MIX}(\Gamma ,X,\mu )}\subseteq{\rm ERG}(\Gamma ,
X,\mu )$.}
\medskip

Consider now the continuous map
\[
a\in A(\Delta ,Y,\nu )\mapsto b(a)\in A(\Gamma ,X\times Y^N,\mu\times\nu^N),
\]
where $b(a)={\rm CInd}(a_0,b_0)^\Gamma_\Delta (a)$.  Then we have,
by 3.1,
\begin{center}
$a$ is mixing $\Rightarrow b(a)$ is mixing
\end{center}
and, by 3.4,
\begin{center}
$b(a)$ is ergodic $\Rightarrow a$ is ergodic.
\end{center}

Let $C=\{a\in A(\Delta ,Y,\nu ): b(a)\in\overline{{\rm MIX}(\Gamma ,X\times
Y^N,\mu\times\nu^N)}\}$.  Then $C$ is closed and MIX$(\Delta ,Y,\nu )
\subseteq C$.  So $\overline{{\rm MIX}(\Delta ,Y,\nu )}\subseteq C$.
If $a\in C$, then $b(a)\in \overline{{\rm MIX}(\Gamma ,X\times
Y^N,\mu\times\nu^N)}\subseteq{\rm ERG}(\Gamma ,X\times Y^N,\mu\times
\nu^N)$.  So $a\in{\rm ERG}(\Delta ,Y,\nu )$, i.e., $\overline{{\rm MIX}
(\Delta ,Y,\nu )}\subseteq{\rm ERG}(\Delta ,Y,\nu )$, and thus $\Delta$
does not have the HAP.

\medskip
{\bf (C)} Let $\Gamma$ be an infinite countable group and let $b_0$
be a free mixing action of $\Gamma$ on $(X,\mu )$.  It is well-known
that there is a free, mixing action $a_0$ of $\bbZ$ on $(X,\mu )$
such that $E=E^X_\bbZ \subseteq F=E^X_\Gamma$. We will construct
below $a_0$ so that moreover $\varphi_E(\gamma )\rightarrow 0$ as
$\gamma\rightarrow\infty$. Then by 3.3, for every action $a$ of
$\bbZ$ on $(Y,\nu )$, the action ${\rm CInd}(a_0, b_0)^\Gamma_\bbZ
(a)$ is mixing, which produces a large supply of seemingly new free,
mixing actions of $\Gamma$.

\medskip
\noindent{\bf Theorem 3.8}. {\it Let $\Gamma$ be an infinite
countable group and let $b_0$ be a free, measure preserving, mixing
action of $\Gamma$ on $(X,\mu )$.  Then there is a free, measure
preserving, mixing action $a_0$ of $\bbZ$ on $(X,\mu )$ with
$E=E^X_\bbZ\subseteq F=E^X_\Gamma$ and $\varphi_E(\gamma
)\rightarrow 0$ as $\gamma\rightarrow\infty$.}
\medskip

{\bf Proof}.  We will use the following two lemmas.
\medskip

\noindent{\bf Lemma 3.9}. {\it Let $\Gamma$ be an infinite countable
group.  Then there exists a positive, symmetric function $f\in
c_0(\Gamma )\setminus\ell^1 (\Gamma )$ such that whenever
$S\subseteq\Gamma$ and $\sum_{\gamma\in S}f(\gamma ) <\infty$, then
for any $\delta_1,\delta_2\in\Gamma ,\sum_{\gamma\in S}f(\delta_1
\gamma\delta_2)<\infty$ (i.e., the summable ideal associated to $f$
is two-sided invariant).}
\medskip

{\bf Proof}.  Fix a sequence $\{Q_n\}_{n\geq 0}$ of finite,
symmetric subsets of $\Gamma$ with $Q_0=\{1\}, Q_n\subseteq Q_{n+1},
Q_{n+1}\setminus (Q_n)^n\neq\emptyset$ and $\bigcup_nQ_n=\Gamma$.
Let $|\gamma |=\min\{n:\gamma\in (Q_n)^n\}$ (this is motivated by an
idea in Struble \cite{ST}).  Note that $|\gamma |=|\gamma^{-1}|$ and
$|\gamma\delta |\leq |\gamma |+|\delta |$, so $|\gamma\delta |\geq
|\ |\gamma | -|\delta |\ |$, for any $\gamma ,\delta\in\Gamma$.  Put
now $f(\gamma )=\tfrac{1} {|\gamma |+1}$.  Then clearly $f\in
c_0(\Gamma )$ and for every $n$, there is $\gamma \in\Gamma$ with
$|\gamma |=n+1$, so $f\not\in\ell^1(\Gamma )$.  Fix now $S\subseteq
\Gamma$ with $\sum_{\gamma \in S}f(\gamma )<\infty$.  For
$\delta\in\Gamma$, let $|\delta |=c$, and notice that
\begin{align*}
\sum_{\gamma\in S}f(\gamma\delta )&=\sum_{\gamma\in
S}\frac{1}{|\gamma\delta |+1} \leq\sum_{\{\gamma :|\gamma |\leq
c\}}\frac{1}{|\gamma\delta |+1}+ \sum_{\{\gamma\in S :|\gamma
|>c\}}\frac{1}{|\gamma |-c+1}
\\
&\leq\sum_{\{\gamma :|\gamma |\leq c\}}\frac{1}{|\gamma\delta
|+1}+\sum_{\gamma\in S} \frac{c+1}{|\gamma |+1}<\infty .
\end{align*}
Similarly, $\sum_{\gamma\in S}f(\delta\gamma )<\infty$ and we are done.\hfill$\dashv$
\medskip

Consider now the given free, measure preserving, mixing action $b_0$ of $\Gamma$
on $(X,\mu )$ with $F=E^X_\Gamma$ the associated equivalence relation.
\medskip

\noindent{\bf Lemma 3.10}. {\it Let $f$ be as in Lemma 3.9.  Let
$R\subseteq F$ be a finite subequivalence relation with uniformly
bounded size of its equivalence classes, and let $A,B\subseteq X$ be
disjoint with $\mu (A)=\mu (B)$ and such that $A\cup B$ is a section
of $R$ {\rm (}i.e., no two distinct members of $A\cup B$ belong to
the same $R$-class{\rm )}.  Suppose also that $\varphi_R(\gamma
)\leq f(\gamma ),\forall \gamma\in\Gamma$.  Then there exists
$\theta\in [[F]]$ with ${\rm dom}(\theta )=A, {\rm rng}(\theta )=B$
such that \setcounter{equation}{0}
\begin{equation}
\varphi_{R\vee\theta}(\gamma )\leq f(\gamma ),\forall\gamma\in\Gamma
.
\end{equation}}
\indent {\bf Proof}.  Let $T\in [F]$ generate $R$ and let $N$ be
such that $T^N=1$.  Notice that for any $\theta\in [[F]]$ with ${\rm
dom}(\theta )\subseteq A,{\rm rng}(\theta )\subseteq B$ we have
\begin{align*}
\varphi_{R\vee\theta}(\gamma )\leq\varphi_R(\gamma )&+\sum^N_{i,j=1}\mu (\{x:
T^i\theta T^j(x)=\gamma\cdot x\})
\\
&+\sum^N_{i,j=1}\mu (\{x:T^i\theta^{-1}T^j(x)=\gamma\cdot x\}).
\end{align*}
Let $\theta\in [[F]]$ be maximal (under inclusion) with ${\rm
dom}(\theta )\subseteq A,{\rm rng}(\theta )\subseteq B$ satisfying
(1).  We will show that this works, i.e., ${\rm dom}(\theta )=A,{\rm
rng}(\theta )=B$.  Otherwise, $A_1=A\setminus {\rm dom}(\theta
),B_1=B\setminus{\rm dom}(\theta )$ have positive measure.

For $\rho\in [[F]],\gamma\in\Gamma$ let
\begin{align*}
s_\gamma (\rho )&=\sum^N_{i,j=1}\mu (\{x:T^i\rho T^j(x)=\gamma\cdot x\})
\\
&\quad\quad +\sum^N_{i,j=1}\mu (\{x:T^i\rho^{-1}T^j(x)=\gamma\cdot x\}),
\end{align*}
and put
\[
S=\{\gamma :\varphi_R(\gamma )+s_\gamma (\theta )\geq f(\gamma )\}.
\]
Let then $K\subseteq\Gamma$ be finite, symmetric such that $\mu
(\bigcup_i\{x: T^i(x)\not\in K\cdot x\})<\tfrac{\mu (A_1)^2}{4}$,
and put $S'=KSK\cup KS^{-1}K$.

\medskip
{\bf Claim}.  $\sum_{\gamma\in S}f(\gamma )<\infty$.
\medskip

{\bf Proof}.  First notice that
$\sum_{\gamma\in\Gamma}\varphi_R(\gamma )=\sum_{\gamma \in
\Gamma}\mu (\{x:\exists i\leq N(T^i(x)=\gamma\cdot
x)\})\leq\sum_{i\leq N}\sum_{\gamma \in \Gamma}\mu
(\{x:T^i(x)=\gamma\cdot x\})\leq N$, as the sets
$\{x:T^i(x)=\gamma\cdot x\},\gamma\in\Gamma$, are pairwise disjoint
by the freeness of the action. Similarly
$\sum_{\gamma\in\Gamma}s_{\gamma}(\theta )\leq 2N^2$, since the sets
$\{x:T^i \theta T^j(x)=\gamma\cdot x\},\gamma\in\Gamma$, are
pairwise disjoint.  Thus $\sum_{\gamma\in S}f(\gamma )\leq
\sum_{\gamma\in S}\varphi_R(\gamma )+\sum_{\gamma \in S}s_\gamma
(\theta )<\infty$.

\medskip
So by Lemma 3.9, $\sum_{\gamma\in S'} f(\gamma)<\infty$ and hence
$\Gamma\setminus S'$ is infinite. Since the action is mixing, there
is $\gamma_0\in\Gamma\setminus S'$ such that $\mu (\gamma_0\cdot
A_1\cap B_1)\geq(3/4)\mu (A_1)^2$.  Then
\[
(K\gamma_0K\cup K\gamma^{-1}_0K)\cap S=\emptyset .
\]
Put
\[
D=(A_1\cap\gamma^{-1}_0(B_1))\setminus (\bigcup_i\{x:T^i(x)\not\in
K\cdot x\} \cup\gamma^{-1}_0\cdot\bigcup_i\{x:T^i(x)\not\in K\cdot
x\}).
\]
Then $\mu (D)\geq (3/4)\mu (A_1)^2-(2/4)\mu (A_1)^2>0$.  Also for
any $i,j,$ $$\ T^i( \gamma_0|D)T^j(x)\in K\gamma_0K\cdot x,$$if the
left-hand side is defined.  Indeed for such $x, T^j(x)\in D$, so $x
= T^{N-j}T^j(x)\in K\cdot T^j(x)$ and $T^i\gamma_0 T^j(x)\in
K\cdot\gamma_0T^j(x)$.  Thus $T^j(x)\in K\cdot x$ and
$T^i\gamma_0T^j(x) \in K \gamma_0 K\cdot x$.  Similarly
$T^i(\gamma^{-1}_0|\gamma_0(D))T^j(x)\in K \gamma^{-1}_0K\cdot x$,
whenever the left-hand side is defined.  In particular,
$$T^i(\gamma_0|D)T^j(x), T^i(\gamma^{-1}_0|\gamma_0(D))T^j(x)$$are
not in $S\cdot x$, when they are defined.  Take now $D'\subseteq D$
with
\[
0<\mu (D')\leq \min_{\gamma\in K\gamma_0K\cup
K\gamma^{-1}_0K}(f(\gamma )-\varphi_R(\gamma )-s_{\gamma}
(\theta))/(2N^2).
\]
This makes sense as $\gamma\in K\gamma_0K\cup K\gamma^{-1}_0K$
implies $\gamma\not\in S$, so $f(\gamma )-\varphi_R(\gamma
)-s_\gamma (\theta )>0$.  Let $\theta_0 = \gamma_0|D'$ and
$\theta'=\theta \cup\theta_0$. We claim that this contradicts the
maximality of $\theta$.  We need to verify that $\theta'$ satisfies
(1).  Now
\begin{align*}
\varphi_{R\vee\theta'}(\gamma )\leq\varphi_{R\vee\theta}(\gamma
)&+\sum^N_{i,j=1}
\mu (\{x:T^i\theta_0 T^j(x)=\gamma\cdot x\}\\
&+\sum^N_{i,j=1}\mu (\{x:T^i(\theta_0)^{-1}T^j(x)=\gamma\cdot x\}).
\end{align*}
But if $T^i\theta_0 T^j(x)=T^i(\gamma_0 |D')T^j(x)=\gamma\cdot x$,
then $\gamma\in K \gamma_0K$, so if $\{x:T^i\theta_0
T^j(x)=\gamma\cdot x\}$ is not empty, then $\gamma \in K\gamma_0K$
and similarly if $\{x:T^i(\theta_0)^{-1}T^j(x)=\gamma\cdot x\}$ is
not empty, $\gamma\in K\gamma^{-1}_0K$.  Thus for any $\gamma\not\in
K\gamma_0K\cup K\gamma^{-1}_0K, \varphi_{R\vee\theta'}(\gamma
)\leq\varphi_{R\vee\theta}(\gamma ) \leq f(\gamma )$ and for
$\gamma\in K\gamma_0K\cup K\gamma^{-1}_0K$,
\begin{align*}
\varphi_{R\vee\theta'}(\gamma ) &\leq\varphi_R(\gamma )+s_\gamma
(\theta )
+2\sum^N_{i,j=1}\mu (D')\\
&\leq\varphi_R(\gamma )+s_\gamma (\theta )+2N^2\frac{f(\gamma
)-\varphi_R
(\gamma )-s_\gamma (\theta )}{2N^2}\\
&\leq f(\gamma ),
\end{align*}
so the proof of the lemma is complete.\hfill$\dashv$

\medskip
To complete the proof of the theorem, we follow the argument of
\cite{Z}, 9.3.2 or \cite{KM}, 7.13, which shows how to construct an
ergodic, hyperfinite subequivalence relation $E\subseteq F$.  That
proof proceeds by constructing a sequence of finite equivalence
relations $E_1\subseteq E_2\subseteq\dots$ each with classes of
bounded size such that $E_1={\rm equality}$ and
$E_{n+1}=E_n\vee\theta_n$, where $\theta_n$ is {\it any} $\theta\in
[[F]]$, with ${\rm dom}(\theta )=A_n, {\rm rng}(\theta )=B_n$, where
$A_n, B_n$ are some appropriately chosen Borel sets with $\mu
(A_n)=\mu (B_n)$ and $A_n\cup B_n$ a section of $E_n$.  Choose now
$f$ as in Lemma 3.9 and by Lemma 3.10 choose inductively $\theta_n$
so that $\varphi_{E_{n+1}}(\gamma )=\varphi_{E_n\vee
\theta_n}(\gamma )\leq f(\gamma ),\forall\gamma$ (we can of course
assume that $f(1) =1$, so that $\varphi_{E_1}(\gamma )\leq f(\gamma
),\forall\gamma$). Then $\varphi_E (\gamma
)=\lim_{n\rightarrow\infty}\varphi_{E_n}(\gamma )\leq f(\gamma )$,
and thus $\varphi_E(\gamma )\rightarrow 0$ as
$\gamma\rightarrow\infty$.

Finally note that, since by Dye's Theorem all ergodic, hyperfinite
equivalence relations are isomorphic, we can find a free mixing
action $a_0$ of $\bbZ$ that induces $E$. \hfill$\dashv$

\medskip
It is an old problem of Schmidt \cite{S} to find out whether there
exist infinite countable groups $\Gamma$ for which every ergodic
action is mixing.  One possible approach towards showing the
non-existence of such groups is the following.  Fix a free, mixing
action $b_0$ of a countable infinite group $\Gamma$ on $(X,\mu )$.
Then there is a free, mixing action $a_0$ of $\bbZ$ on $(X,\mu )$
such that $E=E^X_\bbZ\subseteq F=E^X_\Gamma$. There is a weakly
mixing action $a$ of $\bbZ$ on $(Y,\nu )$ which is not mixing. One
might hope that by constructing judiciously $a_0, {\rm
CInd}(a_0,b_0)^\Gamma_\bbZ(a)$ might also be weakly mixing but not
mixing.

Consider now the case of non-amenable groups $\Gamma$.  Gaboriau and
Lyons \cite{GL} have shown that if $\Gamma$ is a non-amenable group,
there is a free, measure preserving, mixing action $b_0$ of $\Gamma$
on $(X,\mu )$ and a free, measure preserving, ergodic action of
$F_2=\langle a,b\rangle$ on $(X,\mu )$ such that
$E=E^X_{F_2}\subseteq F=E^X_\Gamma$. We will show below that one can
find such a pair of actions so that moreover $\varphi_E(\gamma
)\rightarrow 0$ as $\gamma\rightarrow\infty$.
\medskip

\noindent{\bf Theorem 3.11 (with I. Epstein)}. {\it  Let $\Gamma$ be
a non-amenable countable group.  Then there is a free, measure
preserving, mixing action $b_0$ of $\Gamma$ on $(X,\mu )$ and a
free, measure preserving, ergodic action $a_0$ of $F_2$ as $(X,\mu
)$ with $E=E^X_{F_2}\subseteq F= E^X_\Gamma$ and $\varphi_E(\gamma
)\rightarrow 0$ as $\gamma\rightarrow \infty$.}

\medskip
We will postpone the proof of this theorem until Section 4, {\bf
(D)}, as it will require some ideas from percolation on Cayley
graphs.

\medskip
{\bf (D)} We will finally show how the combination of 3.11 and the
work of Epstein \cite{E}, who showed that any non-amenable countable
group $\Gamma$ has uncountably many non-orbit equivalent free,
measure preserving, ergodic actions, provides the following
strengthening to a non-classification result and also sharpens it by
restricting to mixing actions.
\medskip

\noindent{\bf Theorem 3.12 (with I. Epstein)}. {\it Let $\Gamma$ be
a non-amenable countable group.  Then $E_0$ can be Borel reduced to
{\rm OE} on the space of free, measure preserving, mixing actions of
$\Gamma$ and {\rm OE} on this space cannot be classified by
countable structures.}
\medskip

{\bf Proof}. We fix a free, measure preserving, mixing action $b_0$
of $\Gamma$ on $(X,\mu )$ and a free, measure preserving, ergodic
action $a_0$ of $F_2$ on $(X,\mu )$ with $E=E^X_{F_2}\subseteq
F=E^X_\Gamma$ such that $\varphi_E(\gamma )\rightarrow 0$.  Then for
any action $a\in A(F_2, Y,\nu )$ we have the co-induced action
$b={\rm CInd}(a_0,b_0)^\Gamma_{F_2}(a)\in A(\Gamma ,Z,\rho )$, where
$Z=X\times Y^N,\rho =\mu\times\nu^N$ and the action $a'\in
A(F_2,Z,\rho )$, which gives a subequivalence relation of that given
by $b$.  The action $b$ is free and by 3.3 it is mixing.  Also $a'$
is free.

From Epstein \cite{E}, and the fact that in our case $b$ is ergodic, we also have
the following additional properties:

(*) For any Borel homomorphism $g:Y\rightarrow\bar Y$ of $a$ to a
{\it free} action $\bar a\in A(F_2,\bar Y,\bar\nu )$, we have,
letting $f:Z\rightarrow Y$ be defined by $f(x,\bar y)=\bar y_0$:
\[
\rho (\{z\in Z:\exists\gamma\neq 1(g\circ f(\gamma\cdot z) = g\circ
f(z))\})=0.
\]

(**) For every $a'$-invariant Borel set $A\subseteq Z$ of positive
measure, if $\rho_A=\frac{\rho |A} {\rho (A)}$, is the normalized
restriction of $\rho$ to $A, f_*\rho_A=\nu$, thus $a$ is a factor of
$(a'|A,\rho_A)$.

Consider now the standard action of ${\rm SL}_2(\bbZ )$ on
$(\bbT^2,\lambda )$ with the usual product measure $\lambda$, and
fixing a copy of $F_2$ of finite index in ${\rm SL}_2(\bbZ )$ let
$\bar a_0$ be the restriction of this action to $F_2$.  We will use
the following basic lemma originally proved in Ioana \cite{I}, when
$F_2\leq\Gamma$, but realized to hold as well in the more general
case stated below in Epstein \cite{E}.

\medskip
\noindent{\bf Lemma 3.13}. {\it Let $\Gamma$ be a countable group
and let $\{b_i\}_{i\in I}$ be an uncountable family of {\rm OE}
free, measure preserving, ergodic actions of $\Gamma$ on $(Z,\rho )$
such that for each $i\in I$ there is a free, measure preserving
action $a'_i$ of $F_2$ on $(Z,\rho )$ with the following two
properties:

(i) $E_{a'_i}\subseteq E_{b_i}$ (where $E_c$ is the equivalence
relation induced by an action $c$),

(ii) There is a Borel homomorphism $f_i:Z\rightarrow\bbT^2$ of
$a'_i$ to $\bar a_0$ such that
\[
\rho (\{z:\exists\gamma\neq 1(f_i(\gamma\cdot z)=f_i(z))\})=0.
\]

Then there is an uncountable $J\subseteq I$ such that given any
$i,j\in J$, there are $a'_i,a'_j$-invariant Borel sets $A_i,A_j$ of
positive measure, so that the actions $a'_i|A_i,a'_j|A_j$ are
isomorphic (with respect to the normalized measures $\rho_{A_i},
\rho_{A_j}$).}

\medskip
Denote by ${\rm Irr}(F_2,\h)$ the Polish space of irreducible
unitary representation of $F_2$ on a separable, infinite-dimensional
Hilbert space $\h$ (see, e.g., \cite{K}, Appendix H).  By a result
of Hjorth \cite{H1} (see also \cite{K}, Appendix H) there is a
conjugacy invariant dense $G_\delta$ set $G(\Gamma ,\h)\subseteq{\rm
Irr}(F_2, \h)$ such that the conjugacy action of the unitary group
$U(\h)$ on $G$ is turbulent.  As a consequence, if $\cong$ denotes
isomorphism between representations, $\cong |G(\Gamma ,\h)$ is not
classifiable by countable structures.

Finally for each unitary representation $\pi$ of $F_2$ on $\h$
denote by $a_\pi$ the corresponding Gaussian action of $F_2$ (on a
space $(\Omega ,\tau )$; see \cite{K}, Appendix E).  It has the
following two properties:

(1) $\pi\cong\rho\Rightarrow a_\pi\cong a_\rho$,

(2) If $\kappa^{a_\pi}_0$ is the Koopman representation on
$L^2_0(\Omega ,\tau )$ associated to $a_\pi$, then $\pi\leq
\kappa^{a_\pi}_0$.

Given now any $\pi\in G(\Gamma ,\h)$, consider the (diagonal)
product action $a(\pi )=\bar a_0\times a_\pi$ on
$(\bbT^2\times\Omega ,\lambda\times\tau )=(Y,\nu )$. Let then $b(\pi
)\in A(\Gamma ,Z,\rho )$ be the co-induced action of $a(\pi )$ and
$a'(\pi )$ the associated $F_2$-action.  Thus $\pi\mapsto b(\pi )$
is a Borel function from $G(\Gamma ,\h)$ into the space of free,
measure preserving, mixing actions of $\Gamma$ on $(Z,\rho )$.  Put
for $\pi ,\rho\in G(\Gamma ,\h)$:
\[
\pi R\rho\Leftrightarrow b(\pi ){\rm OE}b(\rho ).
\]
Then $R$ is an equivalence relation on $G(\Gamma ,H)$ and $\pi\cong\rho\Rightarrow
\pi R\rho$.

\medskip
{\bf Claim}. {\it $R$ has countable index over $\cong$.}

\medskip
Granting this, the proof is completed as follows.  First, to see
that $E_0$ can be Borel reduced to OE on the space of free, measure
preserving, mixing actions of $\Gamma$ it is of course enough to
show that it can be Borel reduced to $R$.  The equivalence relation
$R$ is analytic with meager classes (as each $\cong$-class in
$G(\Gamma ,\h)$ is meager in $G(\Gamma ,\h)$, and every $R$-class
contains only countably many $\cong$-classes), so $R$ is meager and
contains the equivalence relation $\cong$ induced by the conjugacy
action of $U(\h)$ on $G(\Gamma , \h)$ which has dense orbits (being
turbulent). Then $E_0$ is Borel reducible to $R$ by the argument in
Becker-Kechris \cite{BK}, 3.4.5.

To prove non-classification by countable structures, it is of course
enough to show that $R$ has the same property.  If this fails,
towards a contradiction, there is Borel $F:G(\Gamma ,\h)\rightarrow
X_L$, where $X_L$ is the standard Borel space of countable models of
a countable language $L$, such that
\[
\pi R\rho\Leftrightarrow F(\pi )\cong F(\rho ),
\]
so that, in particular,
\[
\pi\cong\rho\Rightarrow F(\pi )\cong F(\rho ).
\]
But then, by turbulence, there is a comeager set $A\subseteq
G(\Gamma ,\h)$ and $M_0\in X_L$ with
\[
F(\pi )\cong M_0,\forall\pi\in A.
\]
Since every $R$-class is meager, there are $R$-inequivalent $\pi ,\rho\in A$, so that
$F(\pi )\not\cong F(\rho )$, a contradiction.
\medskip

{\bf Proof of the claim}.  Assume, towards a contradiction, that
there is an uncountable family $\{\pi_i\}_{i\in I}\subseteq G(\Gamma
,\h)$ of pairwise non-isomorphic representations such that if we put
$b_i=b(\pi_i)$, then $\{b_i\}_{i\in I}$ are OE.  Let
$a'_i=a'(\pi_i)$, so that $E_{a'_i}\subseteq E_{b_i}$.  Moreover if
$f_i:Z\rightarrow \bbT^2$ is given by $f_i=g\circ f$, where
$g:Y\rightarrow\bbT^2$ is the projection, then $f_i$ is a Borel
homomorphism of $a'_i$ to $\bar a_0$ such that
\[
\rho (\{z:\exists\gamma\neq 1(f_i(\gamma\cdot z)=f_i(z))\})=0
\]
(by property (*) of the co-induced action mentioned earlier).  So,
by Lemma 3.13, there is an uncountable $J\subseteq I$ such that
given any $i,j\in J$, there are $a'_i,a'_j$-invariant Borel sets
$A_i,A_j$ of positive measure, so that the actions
$a'_i|A_i,a'_j|A_j$ are isomorphic.  Note that we also have, by
property (**) of the co-induced action, that $a(\pi_i)=\bar
a_0\times a_{\pi_i}$ is a factor of $a'_i|A_i$.

Fix any $i_0\in J$.  Then for any $j\in J$, fixing $A_{i_0},A_j$ as
above, $\pi_j\leq \kappa^{a_{\pi_j}}_0\leq \kappa^{a_0\times
a_{\pi_j}}_0\leq \kappa^{a'_j|A_j}_0\cong \kappa^{a'_{i_0}
|A_{i_0}}\leq \kappa^{a'_{i_0}}$.  This produces an uncountable
family $\{\pi_j\}_{j \in J}$ of pairwise non-isomorphic irreducible
subrepresentations of $\kappa^{a'_{i_0}}$, which is
impossible.\hfill$\dashv$

\section{Percolation on Cayley graphs of groups}

{\bf (A)} Let $\Gamma$ be a finitely generated group,
$Q=\{\gamma_1,\dots ,\gamma_n\}$ a set of generators for $\Gamma$,
and let $\g_Q=\langle\Gamma ,\bfe_Q\rangle$ be the left Cayley graph
of $\Gamma$, with respect to $Q$, whose set of edges $\bfe_Q$
consists of all $\{\gamma ,\gamma_i\gamma\},1\leq i\leq
n,\gamma\in\Gamma$. (When we deal with Cayley graphs, we always
assume that $1\not\in Q$.) The group $\Gamma$ acts on $\Gamma$ and
thus on the Cayley graph by right translations: $\delta\cdot\{\gamma
,\gamma_i\gamma\}=\{\gamma\delta^{-1}, \gamma_i\gamma\delta^{-1}\}$.
So $\Gamma$ acts also on the space $\Omega_Q=2^{\bfe_Q}$ by shift:
\[
(\gamma\cdot\omega )(\{\delta ,\epsilon\})=\omega (\gamma^{-1}\cdot\{
\delta ,\epsilon\})=\omega (\{\delta\gamma ,\epsilon\gamma\}).
\]
Every $\omega\in\Omega_Q$ can be viewed as the subgraph of $\g_Q$ with
vertices $\Gamma$ and an edge $\{\gamma ,\gamma_i\gamma\}$ belonging to
$\omega$ iff $\omega (\{\gamma ,\gamma_i\gamma\})=1$.  The connected
components of $\omega$ are called the {\it clusters} of $\omega$.

An {\it invariant bond percolation} on this Cayley graph is a
$\Gamma$-invariant probability Borel measure $\bfp$ on $\Omega_Q$.
The percolation $\bfp$ is {\it ergodic} if the action of $\Gamma$ is
ergodic with respect to $\bfp$.

Consider now a free measure preserving action of $\Gamma$ on $(X,\mu
)$, fix Borel subsets $A_1,\dots ,A_n$ of $X$ such that if $\gamma_j
= \gamma_i^{-1}$, then $\gamma_i \cdot A_i = A_j$, and define $\Phi
:X\rightarrow \Omega_Q$ by
\[
\Phi (x)(\{\delta ,\gamma_i\delta\})=1\Leftrightarrow\delta\cdot x\in A_i.
\]
If $\gamma_i\cdot A_i=B_i$, note that
\begin{align*}
\Phi (x)(\{\delta ,\gamma^{-1}_i\delta\})&=\Phi (x)(\{\gamma^{-1}_i\delta ,
\gamma_i(\gamma^{-1}_i\delta )\})=1\\
&\Leftrightarrow\gamma^{-1}_i\delta\cdot x\in A_i\\
&\Leftrightarrow\delta\cdot x\in B_i.
\end{align*}
It is easy to check that
\[
\Phi (\gamma\cdot x)=\gamma\cdot\Phi (x)
\]
and thus if $\bfp =\Phi_*\mu ,\bfp$ is an invariant bond percolation on the
Cayley graph $\g_Q$.

Note that if $E=\gamma_1|A_1\vee\dots\vee\gamma_n|A_n$ is the
subequivalence relation of $F=E^X_\Gamma$ generated by
$\gamma_1|A_1,\dots ,\gamma_n|A_n$ and if $x\in X$ is such that
$\Phi (x)=\omega$, then $\gamma\mapsto\gamma\cdot x$ is a 1-1
correspondence of $\Gamma$ with $\Gamma\cdot x=[x]_F=\{y:yFx\}$,
which sends the clusters of $\omega$ to the $E$-classes contained in
$[x]_F$. So the structure of the $E$-classes in $[x]_F$ is
equivalent to the structure of clusters of $\Phi (x)=\omega$.  The
$E$-class of $x$ corresponds to the cluster of 1 in $\omega$.  We
let
\[
X_\infty =\{x\in X:[x]_E\text{ is infinite}\}.
\]

For $\gamma ,\delta\in\Gamma$, put
\begin{center}
$C_{\gamma ,\delta}=\{\omega\in\Omega_Q:\gamma ,\delta$ are in the same
cluster of $\omega\}$
\end{center}
and (as in Lyons-Schramm \cite{LS}), let
\[
\tau (\gamma ,\delta )=\bfp (C_{\gamma ,\delta}).
\]
Note that by invariance
\[
\tau (\gamma ,\delta )=\tau (\gamma\epsilon ,\delta\epsilon ),\forall
\epsilon\in\Gamma ,
\]
so $\tau (\gamma ,\delta )=\tau (1,\delta\gamma^{-1})$.  If
\[
C_\gamma =C_{1,\gamma}
\]
and
\[
A_\gamma =\{x:\Phi (x)\in C_\gamma\},
\]
then $x\in A_\gamma$ iff $\gamma\cdot x\in [x]_E$ iff $x\in A_E(\gamma )$
and so
\[
\varphi_E(\gamma )=\mu (A_E(\gamma ))=\bfp (C_\gamma )=\tau (1,\gamma ),
\]
where we identify here $\gamma$ with $x\mapsto\gamma\cdot x$.

\medskip
{\bf (B)} Conversely, let $\bfp$ be an invariant bond percolation on
$\g_Q$. The action of $\Gamma$ on $\Omega_Q$ might not be free.  So
fix a free, measure preserving action of $\Gamma$ on $(Y,\nu )$ and
consider the product action of $\Gamma$ on $(X,\mu )=(\Omega_Q\times
Y,\bfp\times\nu )$, which is clearly free.  (If the action of
$\Gamma$ on $\Omega_Q$ is already free, we can simply take $(X,\mu
)=(\Omega_Q,\bfp )$.)  Let
\begin{align*}
C_i&=\{\omega\in\Omega_Q:\omega (1,\gamma_i)=1\},\\
A_i&=C_i\times Y\subseteq X.
\end{align*}
Consider $\gamma_1|A_1,\dots ,\gamma_n|A_n$ and $\Phi :X\rightarrow\Omega_Q$
defined as before.  Then if $x=(\omega ,y)$, we have
\begin{align*}
\Phi (x)(\{\delta ,\gamma_i\delta\})=1&\Leftrightarrow \delta\cdot x\in A_i\\
&\Leftrightarrow (\delta\cdot\omega ,\delta\cdot y)\in A_i\\
&\Leftrightarrow \delta\cdot\omega\in C_i\\
&\Leftrightarrow \delta\cdot\omega (\{1,\gamma_i\})=1\\
&\Leftrightarrow \omega (\{\delta ,\gamma_i\delta\})=1,
\end{align*}
i.e., $\Phi (x)=\Phi (\omega ,y)=\omega$.

\medskip
{\bf (C)} For further reference, we discuss some additional concepts
and results concerning percolation.

In the context of {\bf (B)}, and assuming that the $\Gamma$-action
on $(X, \mu)$ is ergodic and $\mu(X_{\infty}) >0$, we say that
$\bfp$ has {\it indistinguishable infinite clusters} if $E|X_\infty$
is ergodic (this is not the standard definition but is justified by
Gaboriau-Lyons \cite{GL}, Prop. 5).

For any invariant bond percolation $\bfp$ on $\g_Q$ and edge
$e\in\bfe_Q$, define $\pi_e:\Omega_Q\rightarrow\Omega_Q$ by
$\pi_E(\omega )(e)=\omega\cup\{e\}$. Then we say that $\bfp$ is {\it
insertion-tolerant} if $\bfp (A)>0\Rightarrow\bfp
(\pi_E(A))>0,\forall e \in\bfe_Q$.  An example of
insertion-tolerant, ergodic percolation is Bernoulli percolation.

It is well-known, see Newman-Schulman \cite{NS} or Lyons-Schramm
\cite{LS}, 3.8, that if $\bfp$ is ergodic and insertion-tolerant,
then exactly one of the following happens: $\omega$ has no infinite
clusters, $\bfp$-a.s.; $\omega$ has infinitely many infinite
clusters, $\bfp$-a.s.; $\omega$ has exactly one infinite cluster
$\bfp$-a.s. Thus in the context of {\bf (B)}, if $\bfp$ is
insertion-tolerant, then either the $E$-classes are finite,
$\mu$-a.e., or there are infinitely many infinite $E$-classes in
each $F$-class, $\mu$-a.e., or there is exactly one infinite
$E$-class in each $F$-class, $\mu$-a.e.

Moreover, again in the context of {\bf (B)}, Lyons-Schramm \cite{LS}
and Gaboriau-Lyons \cite{GL}, Prop. 6, show that if $\bfp$ is
ergodic and insertion-tolerant, and has infinite clusters,
$\bfp$-a.s., then $\bfp$ has indistinguishable infinite clusters.

Finally, Lyons-Schramm \cite{LS}, Theorem 4.1, show that for any
ergodic, insertion-tolerant, invariant bond percolation $\bfp$,
\[
\inf_{\gamma\in\Gamma} \tau (1,\gamma )>0\Rightarrow\omega\text{ has
a unique infinite cluster, }\bfp{\rm -a.s.}
\]
This implies that in the context of {\bf (B)}, if $\bfp$ is ergodic
and insertion-tolerant, then if
$\inf_{\gamma\in\Gamma}\varphi_E(\gamma )>0$, there is a unique
infinite $E$-class in each $F$-class, $\mu$-a.e., and thus
$E|X_\infty =F|X_\infty ,\mu$-a.e.

We note here that one can give an alternative proof of Theorem 4.1
in \cite{LS} by using the results in the present paper and the
indistinguishability of infinite clusters for insertion-tolerant
percolations. Indeed assume that ${\rm inf}_{\gamma\in\Gamma}
\tau(1,\gamma)
> 0$. In the notation of {\bf (A), (B)} above, since $\tau (1,\gamma
) = \varphi_E (\gamma)$, this means that $\varphi^0_E ={\rm
inf}_{\gamma\in\Gamma} \varphi_E (\gamma) >0$, so by 2.5, there is
an $E$-invariant Borel set $A\subseteq X$ of positive measure such
that $[F|A:E|A] <\infty$. By taking in {\bf (B)} the action of
$\Gamma$ on $(Y,\nu)$ to be weakly mixing, we have that
$F=E^X_{\Gamma}$ is ergodic and thus $A$ meets every $F$-class. It
follows that $A\subseteq X_{\infty}$ and thus, since $\bfp$ is
insertion-tolerant, so that $E|X_{\infty}$ is ergodic, we have that
$A=X_{\infty}$. So $\omega$ has finitely many infinite clusters,
$\bfp$-a.s. and therefore exactly one.

\medskip
{\bf Remark 4.1}.  Note that for any free, measure preserving action
of an infinite group $\Gamma$ on $(X,\mu )$ and any subequivalence
relation $E\subseteq F=E^X_\Gamma$, if $E$ has finite classes, then
$\varphi_E(\gamma )\rightarrow 0$, as $\gamma \rightarrow\infty$.
Indeed, $\varphi_E(\gamma )=\int f(\gamma , x)d\mu (x)$, where
$f(\gamma ,x)=1$ if $\gamma\cdot x\in [x]_E;=0$, if $\gamma\cdot
x\not\in [x]_E$. Since for each $x, f(\gamma ,x)\rightarrow 0$, this
conclusion follows by the Lebesgue Dominated Convergence Theorem.

\medskip
{\bf (D)} We finish this section by now giving the promised
\medskip

{\bf Proof of Theorem 3.11}. We start with the following lemma.
\medskip

\noindent{\bf Lemma 4.2}.  {\it Consider a free, measure preserving,
ergodic action of a countable group $\Gamma$ on $(X,\mu )$, a Borel
set $A\subseteq X$ of positive measure and $E_A$ a Borel equivalence
relation on $A$ satisfying $E_A\subseteq E^X_\Gamma$ and
$\lim_{\gamma\rightarrow\infty} \mu (\{x\in A:\gamma\cdot x\in A \
\&\ \gamma\cdot xE_Ax\})=0$.  Then there exists a Borel equivalence
relation $E\subseteq E^X_\Gamma$ on $X$ with $E|A=E_A$ such that
$\varphi_E(\gamma ) \rightarrow 0$ as $\gamma\rightarrow\infty$.
Moreover, if $E_A$ is treeable {\rm (}resp., ergodic{\rm)}, then $E$
is treeable {\rm (}resp., ergodic{\rm )}.}
\medskip

{\bf Proof}.  Let $X\setminus
A=\bigsqcup_{\gamma\in\Gamma}D_\gamma$, with $\gamma (D_\gamma
)\subseteq A,\forall\gamma\in\Gamma$.  This exists since $A$ is
complete section for $E^X_\Gamma$.  Let $E=E_A\vee\{\gamma |D_\gamma
:\gamma \in\Gamma\}$ be the equivalence relation generated by $E_A$
and $\{\gamma |D_\gamma :\gamma\in\Gamma\}$.  Clearly $E$ has all
the required properties except perhaps that $\varphi_E(\gamma
)\rightarrow 0$ as $\gamma\rightarrow\infty$, which we now proceed
to verify.

Define $f:X\rightarrow\Gamma$ by
\[
f(x)=\gamma\Leftrightarrow x \in D_\gamma\text{ or }(x\in A\ \&\
\gamma =1).
\]
Fix $\epsilon >0$.  Let $K\subseteq\Gamma$ be finite, symmetric such
that $\mu (\{x:f(x)\not\in K\})<\epsilon$.  Let $F\subseteq\Gamma$
be finite such that $\mu (\{x\in A:\gamma\cdot x\in A\ \&\
\gamma\cdot x E_A x\})<\epsilon /|K|^2$ for $\gamma\not\in F$ (where
$|K|={\rm card}(K)$).  We will show that $\varphi_E (\gamma
)<3\epsilon$ if $\gamma\not\in KFK$.  Indeed, fix such $\gamma$.
Notice that for each $x$ such that $\gamma\cdot xEx,\gamma$ can be
uniquely written as $f(\gamma\cdot x)^{-1}\gamma'_xf(x)$ for some
$\gamma'_x$.  We have
\[
\mu (\{x:\gamma\cdot xEx\})<2\epsilon +\mu (\{x:\gamma\cdot xEx\ \&\ f(x)\in K\
\&\ f(\gamma\cdot x)^{-1}\in K\}).
\]
Now it only remains to notice that if $x$ is in the second set
above, $\gamma'_x \not\in F$ and $\gamma'_x\in K\gamma K$.  So, by
the choice of $F$, the second summand is bounded by $\epsilon$ and
we are done.\hfill$\dashv$
\medskip

By Gaboriau-Lyons \cite{GL}, we fix a free, measure preserving,
mixing action $\bar b_0$ of $\Gamma$ on $(Y,\rho )$ and a free,
measure preserving, ergodic action $\bar a_0$ of $F_2$ on $(Y,\rho
)$ with $E_{\bar a_0}\subseteq E_{\bar b_0}$ and such that moreover
if $F_2=\langle g_1,g_2\rangle$ then $\langle g_1\rangle$ also acts
ergodically on $(Y,\rho )$ and thus, by Dye's Theorem, we can in
fact assume that the action of $\langle g_1\rangle$ is mixing.  Let
$N=[E_{\bar b_0}:E_{\bar a_0}]$.  We can of course assume that
$N>1$.

Let $\g =\langle F_2,\bfe\rangle$ be the Cayley graph of $F_2$
corresponding to the set of generators $Q=\{g_1,g_2\}$.  We consider
on $\g$ Bernoulli $p$-percolation $\bfp_p$ with $1/3 <p<1$, so that
$\omega$ has infinitely many infinite clusters $\bfp_p$-a.s., see
\cite{LP}.  Consider the usual shift action $a$ of $F_2$ on $(
\Omega ,\nu )$, where $\Omega = 2^\bfe $, $\nu =\bfp_p$, and the
co-induced action $b_0= {\rm CInd}(\bar a_0,\bar
b_0)^\Gamma_{F_2}(a)$ on $(X,\mu
)=(Y\times\Omega^N,\rho\times\nu^N)$. This action is mixing by 3.1.
Consider also the associated action $a'$ of $F_2$ on $(X,\mu )$, so
that $E_{a'}\subseteq E_{b_0}=E^X_\Gamma =F$ and $a$ is a factor of
$a'$ via $\varphi (y,\bar\omega )=\bar\omega_0$, so, in particular,
$a'$ is free. The same argument as in the proof of 3.1 also shows
that the action of $\langle g_1\rangle$ on $X$ is mixing, hence $a'$
is ergodic as well.

Recall that the action $a'$ of $F_2$ on $(X,\mu )$ is given (in the notation of
Section 3) by
\[
\delta\cdot (y,\bar\omega )=(\delta\cdot y,(n\mapsto\bar\delta
(\delta, y)_n\cdot\bar\omega_{\pi (\delta,y)^{-1}(n)}),
\]
where $\pi (\delta,y)(k)=n\Leftrightarrow [C_k(y)]_{E_{\bar
a_0}}=[C_n (\delta\cdot y)]_{E_{\bar a_0}}$ and $\bar\delta
(\delta,y)_n\cdot C_{\pi (\delta,y)^{-1}(n)}(y)=C_n(\delta\cdot y)$.
Since $C_0(y)=y$, we have $\pi (\delta,y)(0)=0$.  Moreover
$\bar\delta (\delta,y)_0\cdot y=\delta \cdot y$, so $\bar\delta
(\delta,y)_0=\delta$. It follows that
\[
\delta\cdot (y,\bar\omega )=(\delta\cdot
y,\delta\cdot\bar\omega_0,(n>0\mapsto \bar\delta
(\delta,y)_n\cdot\bar\omega_{\pi (\delta,y)^{-1} (n)})),
\]
so (up to an obvious isomorphism) $X=\Omega\times
(Y\times\Omega^{N\setminus\{0\}})$ and the action $a'$ is the
product action of $a$ and a free, measure preserving action of $F_2$
on
$(Y\times\Omega^{N\setminus\{0\}},\rho\times\nu^{N\setminus\{0\}})$.
The projection of $X=Y\times\Omega^N$ to the first factor $\Omega$
in the above product is of course $\varphi (y,\bar\omega
)=\bar\omega_0$.  Thus we are in the situation of Section 4, {\bf
(B)}.

We can then define the associated subequivalence relation $E'\subseteq E_{a'}$ by
\[
x_1E'x_2\Leftrightarrow\exists\delta\in F_2(\delta\cdot x_1=x_2\ \&\ 1,\delta
\text{ are in the same }\varphi(x_1)\text{ cluster),}
\]
i.e., $E'=g_1|A_{g_1}\vee g_2|A_{g_2}$, where $A_{g_i}=\{x:\varphi
(x)(\{1,g_i\})=1\}$. Let
\begin{align*}
A=X_\infty &=\{x:[x]_{E'}\text{ is infinite\}}\\
&=\{x:\text{ the cluster of 1 in }\varphi (x)\text{ is infinite\}}.
\end{align*}
Since $\bfp_p$ has indistinguishable infinite clusters, $E_A=E'|A$
is ergodic. It is also clear that $E_A$ is treeable and
non-hyperfinite, as the canonical treeing on it has infinitely many
ends, a.e. (see \cite{LP}, 7.29).

We next show that $E_A$ satisfies the hypothesis of 4.2.  Indeed we have by Fubini
\begin{align*}
\mu (\{(y,\bar\omega ): &\gamma\cdot (y,\bar\omega )E_A(y,\bar\omega )\})=\\
&\int\nu^N(\{\bar\omega :\gamma\cdot (y,\bar\omega )E_A(y,\bar\omega
)\})d\rho (y)
\end{align*}
and it is enough to show that the function under the integral
converges to 0 as $\gamma\rightarrow\infty$ for any $y\in Y$.  Fix
$y\in Y$ and consider an arbitrary sequence
$\gamma_n\rightarrow\infty$.  Notice that if $\gamma\cdot
(y,\bar\omega ) E_A(y,\bar\omega )$, then $\gamma\cdot (y,\bar\omega
)E'(y,\bar\omega )$, so there is $\delta\in F_2$ with $\gamma\cdot
(y,\bar\omega )=\delta\cdot (y,\bar\omega )$ and thus $\gamma\cdot
\bar y=\delta\cdot\bar y$.  Since we can clearly assume that for
each $n$, $\{\bar\omega : \gamma_n\cdot (y,\bar\omega
)E_A(y,\bar\omega )\}\neq\emptyset$, there is (unique) $\delta_n\in
F_2$ with $\gamma_n\cdot y=\delta_n\cdot y$ and $\gamma_n\cdot
(y,\bar\omega )E_A(y,\bar\omega )\Rightarrow\gamma_n\cdot (y,\bar
\omega )=\delta_n\cdot (y,\bar\omega )$. Clearly
$\delta_n\rightarrow \infty$. So
\begin{align*}
&\nu^N(\{\bar\omega :\gamma_n\cdot (y,\bar\omega )E_A(y,\bar\omega )\})=\\
&\nu^N(\{\bar\omega :\delta_n\cdot (y,\bar\omega )E_A(y,\bar\omega )\})\leq\tau
(1,\delta_n)\rightarrow 0.
\end{align*}

Thus by 4.2 there is a treeable, ergodic, non-hyperfinite
equivalence relation $E_1\subseteq E^X_\Gamma$ with $E_1|A=E_A$ such
that $\varphi_{E_1}(\gamma )\rightarrow 0$ as
$\gamma\rightarrow\infty$.  In particular, the cost of $E_1$ is
$>1$.  Then by \cite{GL}, Proposition 13, one can find a free,
measure preserving, ergodic action $a_0$ of $F_2$ on $(X,\mu )$ such
that, letting $E=E^X_{F_2}$ be its associated equivalence relation,
we have $E\subseteq E_1\subseteq F=E^X_\Gamma$.  Clearly
$\varphi_E(\gamma )\rightarrow 0$ as $\gamma\rightarrow\infty$ and
the proof of 3.11 is complete.

\section {Property (T) groups}

{\bf (A)} Let now $\Gamma$ be an infinite group with Kazhdan's property (T).  A
{\it Kazhdan pair} for $\Gamma$ is a pair $(Q,\epsilon )$, where $Q$ is a finite
generating set for $\Gamma$ and $\epsilon >0$ is such that for any unitary
representation $\pi$ of $\Gamma$ on a Hilbert space $\h$, if there is a vector
$\xi\in H$ with $\Vert\pi (\gamma )(\xi )-\xi\Vert <\epsilon\Vert\xi\Vert ,\forall
\gamma\in Q$, then there is a non-0 $\Gamma$-invariant vector.  The group $\Gamma$
having property (T) is equivalent to the assertion that there is a Kazhdan pair
$(Q,\epsilon )$ and also equivalent to the assertion that for every finite
generating set $Q$ there is $\epsilon >0$ with $(Q,\epsilon )$ a Kazhdan pair. Let
\[
\epsilon_Q={\rm max}\{\epsilon :(Q,\epsilon )\text{ is a Kazhdan
pair}\}>0
\]
be the {\it maximum Kazhdan constant} associated to $Q$ (this is
sometimes denoted by $\k (Q,\Gamma )$).  It is easy to see that
$\epsilon_Q\leq\sqrt{2}$ and Shalom (private communication) has
shown that $\sup_Q\epsilon_Q=\sqrt{2}$ (where the sup is over all
finite generating sets). We state below a more precise quantitative
version.

\medskip
\noindent{\bf Proposition 5.1 (Shalom)}. {\it Let $\Gamma$ be a
countable group satisfying property {\rm (T)}. Let $(Q,\epsilon )$
be a Kazhdan pair for $\Gamma$, with $Q$ symmetric containing 1. Let
${\rm card}(Q)=k$.  Then for every $n\geq 1$,
\[
\bigg( Q^n,\sqrt{2\biggl[ 1- \big( \frac{k-\epsilon^2/2}{k} \big)^n
\biggr]} \bigg)
\]
is also a Kazhdan pair.  In particular, $\sup_Q\epsilon_Q=\sqrt{2}$,
where the $\sup$ is over all the finite generating sets.}

\medskip
{\bf Proof}.  Note that $(\bar Q,\bar \epsilon )$ is a Kazhdan pair
iff for any unitary representation $\pi :\Gamma\rightarrow U(\h )$
which has no non-0 invariant vectors, we have
\[
\max_{\gamma\in Q}\Vert\pi (\gamma )(\xi )-\xi\Vert\geq\bar\epsilon
,
\]
for every unit vector $\xi\in\h$ or equivalently
\begin{gather}
\min_{\gamma\in Q}{\rm Re}\;\langle\pi (\gamma )(\xi
),\xi\rangle\leq 1-\frac{\bar\epsilon^2}{2}, \tag{$*$}
\end{gather}
for any unit vector $\xi\in\h$.

So fix $Q'=Q^n,\epsilon'=\sqrt{2[1-(\tfrac{k-\epsilon^2/2}{k})^n]}$
in order to show that $(*)$ holds for $(Q',\epsilon')$ and any $\pi$
without invariant non-0 vectors.  For this define the averaging
operator
\[
T=\frac{1}{k}\sum_{\gamma\in Q}\pi (\gamma ).
\]
Since $T$ is self-adjoint, we have
\begin{align*}
\Vert T\Vert &=\sup_{\Vert\xi\Vert =1}|\langle T(\xi ),\xi\rangle |\\
&=\sup_{\Vert\xi\Vert =1}\frac{1}{k}\left |1+\sum_{\gamma\in
Q\setminus\{1\}}{\rm Re}\;
\langle\pi (\gamma )(\xi ),\xi\rangle\right |\\
&\leq\sup_{\Vert\xi\Vert =1}\frac{1}{k}
\left ((k-1)+(1-\frac{\epsilon^2}{2})\right )\\
&=\frac{k-\frac{\epsilon^2}{2}}{k}
\end{align*}
(by $(*)$ for $(Q,\epsilon )$).

Then for every $n\geq 1$ and every unit vector $\xi\in\h$, we have
\[
\frac{1}{k^n}\sum_{\gamma_1,\dots ,\gamma_n\in Q}{\rm
Re}\;\langle\pi (\gamma_1\cdots\gamma_n) (\xi ),\xi\rangle =\langle
T^n(\xi ),\xi\rangle
\]
\[
\leq\Vert T^n\Vert\leq\Vert T\Vert^n\leq \left
(\frac{k-\frac{\epsilon^2}{2}} {k} \right )^n.
\]
This gives
\[
\min_{\gamma\in Q^n}{\rm Re}\;\langle\pi (\gamma )(\xi
),\xi\rangle\leq\left ( \frac{k-\frac{\epsilon^2}{2}} {k} \right
)^n,
\]
for every unit vector $\xi\in\h$.\hfill$\dashv$
\medskip

We can also see that restricting the number of generators to a fixed size $n$ gives
an upper bound strictly less than $\sqrt 2$.  More precisely we have for each infinite
group $\Gamma$ with property (T):
\[
\epsilon_n(\Gamma )=\sup\{\epsilon_Q:Q\text{ generates }\Gamma ,{\rm
card}(Q)\leq n\}\leq\sqrt{2}\cdot\sqrt{\frac{2n-1}{2n+1}}.
\]
To see this fix $Q$ with ${\rm card}(Q)\leq n$ and consider the left
regular representation $\lambda$ of $\Gamma$, let $\bar
Q=Q\cup\{1\}\cup Q^{-1}$, so that ${\rm card}(\bar Q)=m\leq 2n+1$,
and let $\xi\in\ell^2(\Gamma )$ be the unit vector $\xi
=\tfrac{1}{\sqrt m}\chi_{\bar Q}$.  Then, as for each $\gamma\in Q,
\gamma\bar Q\cap\bar Q$ contains $\{1,\gamma\}$, we have
$\langle\gamma\cdot\xi , \xi\rangle =\tfrac{1}{m}{\rm
card}(\gamma\bar Q\cap\bar Q)\geq\tfrac{2}{m}\geq \tfrac{2}{2n+1}$,
so $\Vert\gamma\cdot\xi -\xi\Vert^2=2(1-\langle\gamma\cdot \xi ,\xi
\rangle )\leq 2(1-\tfrac{2}{2n+1})=2(\tfrac{2n-1}{2n+1})$.  Since
$\lambda$ has no non-0 invariant vectors,
$\epsilon_Q\leq\sqrt{2}\cdot\sqrt{\tfrac{2n-1}{2n+1}}$. We do not
know if the above upper bound for $\epsilon_n (\Gamma)$ is best
possible.

We see from 5.1 that if $Q$ is a symmetric generating set with ${\rm
card} (Q) =k$ and $(Q,\epsilon )$ is a Kazhdan pair for $\Gamma$,
then
\[
\epsilon_{k^n}(\Gamma
)\geq\sqrt{2}\cdot\sqrt{1-(\frac{k-\epsilon^2/2}{k})^n}, \forall
n\geq 1.
\]
In a preliminary version of this paper, we raised the following
question: Fix $n\geq 2$. Is $\inf\{\epsilon_n(\Gamma ): \Gamma$ is
an $n$-generated infinite group with property (T)\}$>0$? Shalom
pointed out that this is false, with the counterexamples being
products of two groups with one factor a finite cyclic group. The
question remains open whether there are such counterexamples if one
considers only property (T) groups which are perfect (i.e., groups
equal to their commutator subgroups).

\medskip
The following fact is well-known.
\medskip

\noindent{\bf Proposition 5.2 (see \cite[1.1.8]{BHV})}. {\it Let
$(Q,\epsilon )$ be a Kazhdan pair for $\Gamma$.  Then for any
$\delta
>0$, any unitary representation $(\pi ,\h )$ of $\Gamma$ and
$\xi\in\h$, if $\forall\gamma\in Q (\Vert\pi (\gamma )(\xi
)-\xi\Vert <\delta\epsilon\Vert\xi\Vert )$, then there is a
$\Gamma$-invariant vector $\eta$ with $\Vert\xi -\eta\Vert\leq\delta
\Vert \xi\Vert$.}
\medskip

The next result is a consequence of Deutch-Robertson \cite{DR}, but we will give the
short proof for the reader's convenience.
\medskip

\noindent{\bf Proposition 5.3}. {\it Let $\Gamma$ be a group with property {\rm
(T)}, $(Q,\epsilon )$ a Kazhdan pair for $\Gamma ,\delta >0$, and $\varphi$
a positive-definite function on $\Gamma$ with $\varphi (1)=1$.  Then $\forall
\gamma\in Q({\rm Re}\ \varphi (\gamma )\geq 1-\tfrac{\delta^2\epsilon^2}{2})$
implies $\forall\gamma\in\Gamma ({\rm Re}\ \varphi (\gamma )\geq 1-2\delta^2)$.}
\medskip

{\bf Proof}.  Let $(\pi ,\h ,\xi )=(\pi_\varphi ,\h_\varphi
,\xi_\varphi )$ be the GNS representation of $\Gamma$ associated to
$\varphi$, so that $\langle\pi (\gamma )(\xi ),\xi\rangle =\varphi
(\gamma )$.  In particular $\xi$ is a unit vector as $\varphi
(1)=1$.  Also $\Vert\pi (\gamma )(\xi )-\xi\Vert^2=2(1-{\rm Re}\
\varphi (\gamma ))$.  So if ${\rm Re}\ \varphi (\gamma )\geq
1-\tfrac{\delta^2\epsilon^2}{2}$ for $\gamma\in Q$, then $\Vert\pi
(\gamma )(\xi )-\xi\Vert^2\leq\delta^2\epsilon^2$, so there is a
$\Gamma$-invariant vector $\eta$ with $\Vert\xi -\eta\Vert\leq
\delta$.  Then $\Vert\pi (\gamma )(\xi )-\eta\Vert\leq\delta$, so
$\Vert\pi (\gamma )(\xi )-\xi\Vert\leq 2\delta$, thus $2(1-{\rm Re}\
\varphi (\gamma ))\leq 4\delta^2, \forall\gamma\in\Gamma$ or ${\rm
Re}\ \varphi (\gamma )\geq 1-2\delta^2,\forall\gamma
\in\Gamma$.\hfill$\dashv$
\medskip

{\bf (B)} Consider now a measure preserving action of a group
$\Gamma$ with property (T) on a standard measure space $(X,\mu )$
and let $F=E^X_\Gamma$ be the associated equivalence relation and
$E\subseteq F$ a subequivalence relation.  Applying the preceding to
$\varphi_E$, we obtain:
\medskip

\noindent{\bf Corollary 5.4}. {\it Let $\Gamma$ be an infinite group
with property {\rm (T)} and $(Q,\epsilon )$ a Kazhdan pair for
$\Gamma$.  If $\Gamma$ acts in a measure preserving way on $(X,\mu
)$ with associated equivalence relation $F=E^X_\Gamma$ and
$E\subseteq F$ is a subequivalence relation, then:

(i) For any $\delta >0 ,\min_{\gamma\in Q}\varphi_E(\gamma )\geq
1-\tfrac{\delta^2 \epsilon^2}{2}$, implies
$\varphi^0_E=\inf_{\gamma\in \Gamma}\varphi_E(\gamma )\geq 1-
2\delta^2$.

(ii) If $\min_{\gamma\in Q}\varphi_E(\gamma
)>1-\tfrac{\epsilon^2}{2}$, then there is an $E$-invariant set of
positive measure $A$ such that $[F|A:E|A]< \infty$.  If moreover
$\min_{\gamma\in Q}\varphi_E(\gamma )>1-\tfrac{\epsilon^2} {4}$,
then $\varphi^0_E>0$ and $[F|A:E|A]\leq\tfrac{1}{\varphi^0_E}$.  If
$\min_{\gamma\in Q}\varphi_E(\gamma )>1-\tfrac{\epsilon^2}{8}$, then
$\varphi^0_E
>\tfrac{1}{2}$ and $F|A=E|A$.  Finally if $\min_{\gamma\in Q}\varphi_E(\gamma )>
1-\tfrac{\epsilon^2}{16}$, then $\varphi^0_E>\tfrac{3}{4}$ and we
can find such an $A$ with $\mu(A)\geq 4\varphi^0_E-3$.

(iii) If the action of $\Gamma$ is free and $E$ is induced by a free
action of $\Delta$, then if $\min_{\gamma\in Q}\varphi_E(\gamma
)>1-\tfrac{\epsilon^2}{2},\ \Gamma \ {\rm and}\ \Delta$ are {\rm ME}
and so $\Delta$ has property {\rm (T)}.}

\medskip
{\bf Proof}.  (i) follows from 5.3.  For (ii) first notice that if
$\min_{\gamma\in Q} \varphi_E(\gamma )>1-\tfrac{\epsilon^2}{2}$,
then for the representation $\tau$ discussed in Section 2, {\bf (A)}
and letting $\tau (\gamma )(\xi )=\gamma\cdot \xi$, we have
$\Vert\gamma\cdot\xi_0-\xi_0\Vert^2=2(1-\langle\gamma\cdot\xi_0,
\xi_0\rangle )=2(1-\varphi_E(\gamma ))<\epsilon^2$, for all
$\gamma\in Q$, so $\tau$ has an invariant non-0 vector, thus, by
2.3, there is an $E$-invariant set $A\subseteq X$ of positive
measure for which $[F|A:E|A]<\infty$.  If $\min_{\gamma\in Q}
\varphi_E(\gamma )>1-\tfrac{\epsilon^2}{4}$, then, by 5.3,
$\varphi^0_E>0$, so by 2.5 we can find such an $A$ with
$[F|A:E|A]\leq\tfrac{1}{\varphi^0_E}$.  If $\min_{\gamma\in
Q}\varphi_E(\gamma )>1-\tfrac{\epsilon^2}{8}$, then again by 5.3,
$\varphi^0_E>\tfrac{1}{2}$, so such an $A$ can be found with
$E|A=F|A$. Finally if $\min_{\gamma\in Q}\varphi_E(\gamma
)>1-\tfrac{\epsilon^2}{16}$, then $\varphi^0_E > \tfrac{3}{4}$ and
such an $A$ can be found with $\mu (A)>4\varphi^0_E-3$, using 2.15.
Clearly (iii) follows from the above and 2.8.\hfill$\dashv$
\medskip

We next note the following quantitative version of 3.4 for groups
with property (T).

\medskip
\noindent{\bf Proposition 5.5}. {\it In the notation of Section 2
{\bf (A)}, let $\Gamma$ have property {\rm (T)} and let $(Q,\epsilon
)$ be a Kazhdan pair for $\Gamma$.  If $\min_{\gamma\in
Q}\varphi_E(\gamma )>1-\tfrac{\epsilon^2}{2}$, then for any $a\in
A(\Delta ,Y,\nu )$, if $b={\rm CInd}(a_0,b_0)^\Gamma_\Delta (a)$,
then
\begin{center}
$b$ is ergodic $\Rightarrow a$ is ergodic.
\end{center}}

{\bf Proof}.  Assume that $a$ is not ergodic and repeat the proof of
3.4, with $k=n=0$. Then for $\xi =f^{(0)},\langle\gamma\cdot\xi
,\xi\rangle =\varphi_E(\gamma^{-1})= \varphi_E(\gamma
),\forall\gamma\in Q$.  So $\Vert\gamma\cdot\xi -\xi\Vert^2=2
(1-\langle\gamma\cdot\xi ,\xi\rangle )<\epsilon^2,\forall\gamma\in
Q$, thus there is a non-0 $\Gamma$-invariant vector, so $b$ is not
ergodic.\hfill$\dashv$
\medskip

{\bf (C)} We next consider some consequences concerning percolation
on Cayley graphs of property (T) groups.

\medskip
\noindent{\bf Theorem 5.6}. {\it Let $\Gamma$ be an infinite group
with property {\rm (T)}, $(Q,\epsilon )$ a Kazhdan pair, and $\bfp$
an invariant, ergodic, insertion-tolerant bond percolation on
$\g_Q$. Then if the survival probability ${\bfp} (\{\omega :\omega
(e)=1 )$ of each edge $e$ is $>1-\tfrac{\epsilon^2}{2}, \omega$ has
a unique infinite cluster, ${\bfp}$-a.s.}

\medskip
{\bf Proof}.  In the context of Section 4, whose notation we keep
below, take the free action of $\Gamma$ on $(Y,\nu )$ to be weakly
mixing, so that the $\Gamma$-action on $(X,\mu )=(\Omega_Q\times
Y,{\bfp}\times\nu )$ is free and ergodic.  If $Q=\{\gamma_1, \dots
,\gamma_n\}$, then we have ${\bfp} (\{\omega :\omega
(\{1,\gamma_i\})=1\})\leq\tau (1,\gamma_i)=\varphi_E(\gamma_i)$, so
$\min_{\gamma\in Q}\varphi_E(\gamma )>1- \tfrac{\epsilon^2}{2}$,
thus, by 5.4 (ii), there is an $E$-invariant set $A\subseteq
X_\infty$ of positive measure with $[F|A:E|A]<\infty$.  Since
${\bfp}$ is insertion-tolerant, by Section 4, {\bf (C)},
$E|X_\infty$ is ergodic, so $A=X_\infty$, which again by Section 4,
{\bf (C)} implies that $[F|A:E|A]=1$, i.e., $F|X_\infty
=E|X_\infty$. This means that $\omega$ has a unique infinite
cluster, ${\bfp}$-a.s. \hfill$\dashv$

\medskip
In the case of Bernoulli bond percolation ${\bfp}_p,p\in (0,1)$ on
$\g_Q$, let $p_u=p_u (Q)$ be the critical probability for uniqueness
defined by:  $p_u = \inf \{ p :\text{there is a unique infinite
cluster}, {\bfp}_p{\rm -a.s.}\}$. In Lyons-Schramm \cite{LS} the
authors show that for property (T) groups $\Gamma$ and any finite
set of generators $Q$ one has $p_u(Q)<1$. From the preceding result
one has a quantitative version.

\medskip
\noindent{\bf Corollary 5.7}. {\it Let $\Gamma$ be an infinite group
with property {\rm (T)} and $(Q,\epsilon )$ a Kazhdan pair.  Then
$p_u (Q)\leq 1-\tfrac{\epsilon^2} {2}$.}
\medskip

We also have the following:
\medskip

\noindent{\bf Corollary 5.8}. {\it For each $\rho >0$ and every
infinite group $\Gamma$ with property {\rm (T)}, there is a finite
set of generators $Q$ for $\Gamma$ such that for any invariant,
ergodic, insertion-tolerant bond percolation ${\bfp}$ on $\g_Q$, if
the survival probability of each edge is $\geq\rho$, then $\omega$
has a unique infinite cluster, ${\bfp}$-a.s.}

\medskip
{\bf Proof}.  By 5.1, $\sup_Q\epsilon_Q=\sqrt{2}$, where the sup is
taken over all finite generating sets $Q$ for $\Gamma$.

So fix $\rho >0,\Gamma$ an infinite group with property (T) and $Q$
a finite generating set of $\Gamma$ such that $\epsilon
=\epsilon_Q>\sqrt{2(1-\rho )}$.  Then for ${\bfp}$ as in the
statement of the present corollary, the survival probability of each
edge is bigger than $1-\tfrac{\epsilon^2}{2}$, so $\omega$ has a
unique infinite cluster, ${\bfp}$-a.s., by 5.6.\hfill$\dashv$

\medskip
There is also a version of 5.8 for arbitrary invariant, ergodic bond
percolations.

\medskip
\noindent{\bf Corollary 5.9}. {\it Let $\Gamma$ be an infinite group
with property {\rm (T)} and $(Q,\epsilon )$ a Kazhdan pair.  Then
for any invariant, ergodic bond percolation ${\bfp}$ on $\g_Q$, if
the probability of survival of every edge is $>1-
\tfrac{\epsilon^2}{2}$, there is $n\geq 1$ and a $\Gamma$-invariant
Borel map $\calc_0:\Omega_Q\rightarrow [2^\Gamma ]^n$ (= the space
of $n$-element subsets of the power set of $\Gamma$) such that
$\calc_0(\omega )$ is a set of $n$ infinite clusters of $\omega
,{\bfp}$-a.s.

In particular, for every $\rho >0$ and every infinite group $\Gamma$
with property {\rm (T)}, there is a finite set of generators $Q$ for
$\Gamma$ such that for any invariant, ergodic bond percolation
${\bfp}$ on $\g_Q$, if the survival probability of each edge is
$\geq\rho$, then we can assign in a $\Gamma$-invariant Borel way a
finite set (of fixed size) of infinite clusters to each $\omega
,\bfp$-a.s.}
\medskip

{\bf Proof}.  We follow the proof of 5.6, whose notation and that of
Section 4 we use below. Let $Q=\{\gamma_1,\dots ,\gamma_n\}$. Then
$\mu (A_i)>1-\tfrac{\epsilon^2}{2}$, so if $E$ is the equivalence
relation induced by $\gamma_1|A_1,\dots ,\gamma_n|A_n$, we have
$\varphi_E(\gamma_i)>1-\tfrac{\epsilon^2}{2}$.  So, by 5.4, there is
an $E$-invariant set $A\subseteq X_\infty$ of positive measure such
that $[F|A:E|A]= k<\infty$.  For each $x=(\omega ,y)$, let
$f(x)=\{\gamma :\gamma\cdot x\in A\}\in 2^\Gamma$.  Clearly for
$\delta \in\Gamma,\ f(\delta\cdot x)=f(x)\delta^{-1}=\delta \cdot
f(x)$, where $\Gamma$ acts on the set of subsets $2^\Gamma$ of
$\Gamma$ by right multiplication.  Moreover $f(x)$ is the union of
$k$ infinite clusters of $\omega$.  Thus the map $f_\omega
:Y\rightarrow [\calc (\omega )]^k=$ (the set of $k$-element subsets
of $\calc (\omega ))$, where $\calc (\omega )=$ (the set of infinite
clusters of $\omega$), given by $f_\omega (y)=$ (the set of clusters
contained in $f(\omega ,y)$), induces a measure $(f_\omega
)_*\nu=\nu_\omega$ on $[\calc (\omega )]^k$.  Moreover,
$\gamma\cdot\nu_\omega =\nu_{\gamma\cdot\omega}$, where
$\gamma\cdot\nu_\omega (B)=\nu_\omega (\gamma^{-1}\cdot B)$, for
$B\subseteq [\calc (\omega )]^k$.  But $(f_\omega )_*\nu$ is a
countably additive measure on the countable set $[\calc (\omega
)]^k$, thus can be viewed as given by a weight function $W(\omega
,\bar C),\bar C\in [\calc (\omega )]^k$, where $0\leq W(\omega ,
\bar C)\leq 1$, and $\sum_{\bar C\in [\calc (\omega )]^k}W(\omega
,C)=1$. Moreover, $W(\omega ,\bar C)=W(\gamma\cdot\omega
,\gamma\cdot\bar C)$.  Let $\{\bar C_1 (\omega ),\dots ,\bar
C_{n(\omega )}(\omega)\}$ be the set of $k$-element subsets of
$C(\omega )$ of maximum weight.  Again $\omega\mapsto\{\bar
C_1(\omega ),\dots ,\bar C_{n(\omega )} (\omega )\}$ is
$\Gamma$-invariant.  By the ergodicity of ${\bfp} ,n(\omega )=n_0,
{\bfp}$-a.s.  Let $\calc_0(\omega )=\bar
C_1(\omega)\cup\dots\cup\bar C_{n_0}(\omega )\in [\calc (\omega
)]^{<\bbN}=$ (the set of finite subsets of $\calc (\omega )$). Again
$\calc_0(\omega )$ is $\Gamma$-invariant, so for some $n,
\calc_0(\omega )\in [C(\omega )]^n$, $\bfp$-a.s., and the proof is
complete.

\hfill$\dashv$

\medskip
{\bf (D)} Next we derive some upper bounds for the cost of a group
with property (T). Below $C(\Gamma )$ denotes the cost of a group
$\Gamma$.  If $\Gamma$ is infinite and has property (T) and $Q$ is a
finite set of generators with card$(Q)=n$, then it is well known
that $1\leq C(\Gamma )<n$ (the strict inequality follows from
Gaboriau \cite{G1}, since no free measure preserving action of
$\Gamma$ is treeable, see Adams-Spatzier \cite{AS}). We prove below
some improvements on this upper bound. It should be pointed out
however that at this time no property (T) groups $\Gamma$ with
$C(\Gamma )>1$ are known to exist.

The next result is obtained by a combination of Lyons-Peres-Schramm
\cite {LPS} and 5.7.

\medskip
\noindent{\bf Theorem 5.10}. {\it Let $\Gamma$ be an infinite group
with property {\rm (T)}.  Let $(Q,\epsilon )$ be a Kazhdan pair for
$\Gamma$.  Then if $n={\rm card}(Q)$,
\[
C(\Gamma )\leq n(1-\frac{\epsilon^2}{2})+\frac{n-1}{2n-1}.
\]}
\indent {\bf Proof}.  For the Bernoulli bond percolation $\bfp_p
,p\in (0,1)$ on $\g_Q$, let $p_c=p_c(Q)$ be the critical probability
for infinite clusters defined by: $p_c = \sup \{ p: \text {all
clusters are finite}, \bfp_{p}{\rm -a.s.}\}$. If $d$ is the degree
of $\g_Q$, so that $d\leq 2n$, it is well-known that
$p_c\geq\tfrac{1}{d-1}$, see Lyons-Peres \cite{LP}. Consider the
probability space
\[
\tilde\Omega_Q=\{\tilde\omega\in [0,1]^{\bfe_Q}:\text{ all edge labels are distinct}\}
\]
equipped with the product measure of the Lebesgue measure on [0,1].
Consider the random variable $\Im :\tilde\Omega_Q\rightarrow
2^{\bfe_Q}$ corresponding to the free minimal spanning forest as
defined in Lyons-Peres-Schramm \cite{LPS}.

Then we have
\begin{equation}
C(\Gamma )\leq\frac{1}{2}\bfe ({\rm deg}_1\Im
)\leq\frac{1}{2}(2+d\int^{p_u}_{p_c} \theta (p)^2dp),
\end{equation}
where ${\rm deg}_1\Im$ is the degree of the identity of $\Gamma$ in
the forest. Here $\theta (p)$ is the probability that the cluster of
1 is infinite in Bernoulli $p$-percolation.  The first inequality,
due to Lyons,  follows from the following observation.  Fix a
positive number $\bar\epsilon$ and consider the probability space
$2^{\bfe_Q}$ equipped with the Bernoulli measure
$\bfp_{\bar\epsilon}$. Consider also the diagonal action of $\Gamma$
on $\tilde\Omega_Q\times 2^{\bfe_Q}$ with associated product measure
$\mu$.  Define the graphing $\g$ of the orbit equivalence relation
as follows:
\begin{align*}
\{(\tilde\omega_1,\omega_1),(\tilde\omega_2,\omega_2)\}\in\g\Leftrightarrow
&\exists\gamma\in Q[\gamma\cdot (\tilde\omega_1,\omega_1)=(\tilde\omega_2,\omega_2)\\
&\text{and }(\{1,\gamma\}\in\Im (\tilde\omega_1)\text{ or
}\{1,\gamma\} \in\omega_1)].
\end{align*}
That $\g$ spans the equivalence relation follows from \cite{LPS},
Theorem 3.22 and for the cost, we have
\[
C_\mu(\g)\leq\frac{1}{2}(\bfe({\rm deg}_1\Im (\tilde\omega ))+\bfe
({\rm deg}_1\omega ))=\frac{1}{2}(\bfe ({\rm deg}_1\Im
)+\bar\epsilon d).
\]
Since $\bar\epsilon$ was arbitrary, we obtain the desired
inequality. The second inequality in (2) is \cite{LPS}, Corollary
3.24.

Thus we have, using 5.7,
\begin{align*}
C(\Gamma ) &\leq 1+\frac{d}{2}(p_u-p_c)\\
&\leq 1+\frac{d}{2}\biggl (\left (1-\frac{\epsilon^2}{2}\right )-\frac{1}{d-1}\biggr )\\
&=\frac{d}{2}\left (1-\frac{\epsilon^2}{2}\right )+\frac{d-2}{2(d-1)}\\
&\leq n\left (1-\frac{\epsilon^2}{2}\right )+\left (\frac{n-1}{2n-1}\right ).
\end{align*}
\hfill$\dashv$

\medskip
When $\Gamma$ is torsion-free we also obtain some additional
estimates.
\medskip

\noindent{\bf Theorem 5.11}. {\it Let $\Gamma$ be an infinite group
with property {\rm (T)} and $(Q,\epsilon )$ a Kazhdan pair for
$\Gamma$. Let ${\rm card}(Q)=n$. If $Q$ contains an element of
infinite order, then
\[
C(\Gamma )\leq n-(n-1)\frac{\epsilon^2}{8}.
\]}
\vskip -5pt {\bf Proof}.  Consider a free, ergodic action of
$\Gamma$ on $(X,\mu )$ with associated equivalence relation
$F=E^X_\Gamma$.  Let $Q=\{\gamma_1,\dots ,\gamma_n\}$, where
$\gamma_1$ has infinite order.  Fix $\delta <\tfrac{1}{2}$ and let
$A\subseteq X$ have measure $\mu
(A)=1-\tfrac{\delta^2\epsilon^2}{2}$.  Let $E=\gamma_1\vee\gamma_2
|A\vee\dots \vee\gamma_n|A$ be the equivalence relation generated by
$\gamma_1,\gamma_2|A, \dots ,\gamma_n|A$, so that $E$ is aperiodic.
Then $C_\mu (E)\leq 1+(n-1)(1-\tfrac {\delta^2\epsilon^2}{2})$. Also
$\varphi_E(\gamma )\geq 1-\tfrac{\delta^2\epsilon^2}{2},
\forall\gamma\in Q$, so, by 5.4 (i), $\varphi^0_E\geq
1-2\delta^2>\tfrac{1}{2}$. Then, by 2.13, $C_\mu (F)\leq C_\mu
(E)\leq 1+(n-1)(1-\tfrac{\delta^2\epsilon^2}{2})$. Taking
$\delta\rightarrow\tfrac{1}{2}$ we are done.\hfill$\dashv$
\medskip

\noindent{\bf Theorem 5.12}. {\it Let $\Gamma$ be an infinite group
with property {\rm (T)} and let $(Q,\epsilon )$ be a Kazhdan pair
for $\Gamma$ with $Q$ containing an element of infinite order.  Then
\[
C(\Gamma )\leq n-\frac{\epsilon^2}{2}.
\]}
\indent {\bf Proof}.  Let $b_0$ be a free, mixing action of $\Gamma$
on $(X,\mu )$ and put $F=E^X_\Gamma$.  Let $Q=\{\gamma_1,\dots
,\gamma_n\}$, where $\gamma_1$ has infinite order.  Consider the
graphing of $F$ given by $\gamma_1,\dots ,\gamma_n$.  Applying the
argument in Kechris-Miller \cite{KM} and (independently) Pichot
\cite{P}, we obtain a treeing of a subequivalence relation
$E\subseteq F$, generated by $\gamma_1, \gamma_2|A_2,\dots
,\gamma_n|A_n$, for some Borel sets $A_2,\dots ,A_n$, such that
$C_\mu (E)\geq C_\mu (F)$.  Thus $C_\mu (F)\leq 1+\sum^n_{i=2}\mu
(A_i)$.  Now $E$ is treeable and ergodic, so by Hjorth \cite{H3},
$E$ is induced by a free action $a_0$ of a group $\Delta$.  If we
had $\mu (A_i)>1-\tfrac{\epsilon^2}{2}$, for all $i=2,\dots ,n$,
then $\min_{\gamma\in Q}\varphi_E(\gamma )>1-\tfrac{\epsilon^2}{2}$,
so, by 5.4 (iii), $\Delta$ has property (T) a contradiction.  So for
some $i=2,\dots ,n,\mu (A_i)\leq 1-\tfrac{\epsilon^2}{2}$, therefore
$C_\mu (E)\leq 1+(n-2)+
(1-\tfrac{\epsilon^2}{2})=n-\tfrac{\epsilon^2}{2}$.

\hfill$\dashv$

\medskip
{\bf (E)} Finally, we note that there is an analog of 5.4 (iii),
when $\Gamma$ does not have the HAP.
\medskip

\noindent{\bf Proposition 5.13}. {\it Let $\Gamma$ be an infinite group without the
{\rm HAP}.  Then there exists $\epsilon >0$ and a finite set $Q\subseteq\Gamma$ with
the following property:

Let $\Delta$ be a group and consider two free, measure preserving
actions of $\Gamma$ and $\Delta$ on $(X,\mu )$ such that
$E=E^X_\Delta\subseteq F=E^X_\Gamma$.  If $\min_{\gamma\in Q}
\varphi_E(\gamma )>1-\epsilon$, then $\Delta$ does not have the {\rm
HAP}.}
\medskip

{\bf Proof}.  Since $\Gamma$ does not have the HAP, we can find
$\epsilon >0$ and $Q\subseteq\Gamma$ finite such that if $\varphi
:\Gamma\rightarrow\bbC$ is a positive-definite function with
$\varphi (1)=1$ and $\varphi\in c_0(\Gamma )$, then $\min_{\gamma
\in Q}\varphi (\gamma )\leq 1-\epsilon$.

Now let $\Delta$ be as above and assume that it has the HAP.  Let
$\psi_n:\Delta \rightarrow\bbC$ be positive-definite functions such
that $\lim_{n\rightarrow\infty} \psi_n(\delta )=1$, for all
$\delta\in\Delta$, and $\psi_n\in c_0(\Delta )$, for all $n$.  If
$A_{\gamma ,\delta}=\{x\in X : \gamma\cdot x=\delta\cdot x\}$, then
the formula $\varphi_n(\gamma )=\sum_{\delta\in\Delta}\psi_n(\delta
)\mu (A_{\gamma , \delta}$) defines a sequence of positive-definite
functions on $\Gamma$.

Next we have that
\[
\lim_{n\rightarrow\infty}\varphi_n(\gamma
)=\sum_{\delta\in\Delta}\mu (A_{\gamma ,\delta} )=\mu (\{x\in X :
\gamma\cdot x\in\Delta\cdot x\}=\varphi_E(\gamma ),\forall\gamma\in
\Gamma .
\]
Thus, to get a contradiction to the non-HAP assumption, it suffices
to show that $\varphi_n\in c_0(\Gamma )$, for all $n$.  This is
clear, since for a fixed $n$ we have that
$\lim_{\delta\rightarrow\infty}\psi_n(\delta )=0$,
$\lim_{\gamma\rightarrow \infty}\mu (A_{\gamma ,\delta})=0$, for all
$\delta\in\Delta$, and $\sum_{\delta} \mu (A_{\gamma, \delta}) \leq
1$, for all $\gamma \in \Gamma$.\hfill$\dashv$

\noindent Department of Mathematics

\noindent Caltech

\noindent Pasadena, CA  91125

\medskip
\noindent \texttt{aioana@caltech.edu, kechris@caltech.edu,
todor@caltech.edu}

\end{document}